\documentclass[mnsc,nonblindrev]{informs3_modified1}

\OneAndAHalfSpacedXI 

\usepackage{natbib}
 \bibpunct[, ]{(}{)}{,}{a}{}{,}%
 \def\bibfont{\small}%
\TheoremsNumberedThrough     
\ECRepeatTheorems

\EquationsNumberedThrough    

\MANUSCRIPTNO{MS-0001-1922.65}

\begin{document}


\RUNAUTHOR{Lam}

\RUNTITLE{On the Impossibility of Statistically Improving Empirical Optimization}

\TITLE{On the Impossibility of Statistically Improving Empirical Optimization: A Second-Order Stochastic Dominance Perspective}

\ARTICLEAUTHORS{%
\AUTHOR{Henry Lam}
\AFF{Department of Industrial Engineering and Operations Research, Columbia University, New York, NY 10027, \EMAIL{henry.lam@columbia.edu}} 
} 

\ABSTRACT{%
When the underlying probability distribution in a stochastic optimization is observed only through data, various  data-driven formulations have been studied to obtain approximate optimal solutions. We show that no such formulations can, in a sense, theoretically improve the statistical quality of the solution obtained from empirical optimization. We argue this by proving that the first-order behavior of the optimality gap against the oracle best solution, which includes both the bias and variance, for any data-driven solution is second-order stochastically dominated by empirical optimization, as long as suitable smoothness holds with respect to the underlying distribution. We demonstrate this impossibility of improvement in a range of examples including regularized optimization, distributionally robust optimization, parametric optimization and Bayesian generalizations. We also discuss the connections of our results to semiparametric statistical inference and other perspectives in the data-driven optimization literature.

}%



\maketitle

\section{Introduction}
We consider a stochastic optimization problem in the form:
\begin{equation}
\min_{x\in\mathcal X}\{Z(x):=\psi(x,P)\}\label{main}
\end{equation}
where $x$ is the decision variable in a known feasible region $\mathcal X\subset\mathbb R^d$, and $Z:\mathbb R^d\to\mathbb R$ is the objective function that depends on an underlying probability distribution $P$ and we denote $Z(x)=\psi(x,P)$. A primary example of $\psi(x,P)$ is the expected value objective function $E_P[h(x,\xi)]$ where $E_P[\cdot]$ denotes the expectation with respect to $P$ that generates a random object $\xi\in\Xi$. This, however, can be more general, including for instance the (conditional) value-at-risk of $\xi$.


We focus on \emph{data-driven optimization} where $P$ is not known but only observed via i.i.d. data. In this situation, the decision maker obtains a data-driven solution, say $\hat x$, typically by solving some reformulation of \eqref{main} that utilizes the data. Let $x^*$ be an (unknown) optimal solution for \eqref{main}. We are interested in the statistical properties of the \emph{optimality gap} or regret
\begin{equation}
\mathcal G(\hat x):=Z(\hat x)-Z(x^*)\label{opt gap}
\end{equation}
This captures the suboptimality of $\hat x$ relative to the oracle optimal solution $x^*$. \eqref{opt gap} is a natural metric to evaluate the quality of an obtained solution and, in the language of statistical learning, it measures the generalization performance relative to the oracle best in terms of the objective value. Obviously, if a solution bears a smaller optimality gap than another solution, then its true objective value or generalization performance is also better by the same magnitude. 

Data-driven optimization as discussed above arises ubiquitously across operations research and machine learning, where the objective $\psi(x,P)$ ranges from an expected business revenue to the loss of a statistical model. To obtain $\hat x$ from data, the arguably most straightforward approach is empirical optimization (EO), namely by replacing the unknown true $P$ with the empirical distribution $\hat P$ in the objective $\psi(x,P)$. For example, in expected value optimization where $\psi(x,P)=E_P[h(x,\xi)]$, this corresponds to minimizing $E_{\hat P}[h(x,\xi)]$, also known as the sample average approximation (SAA) \citep{shapiro2014lectures}. Other than EO, there are plenty of actively studied reformulations, including regularization that adds penalty terms to the objective \citep{friedman2001elements}, and data-driven distributionally robust optimization (DRO) \citep{delage2010distributionally,goh2010distributionally,ben2013robust,wiesemann2014distributionally,lim2006model,rahimian2019distributionally} where one turns the objective $\psi(x,P)$ into $\max_{Q\in\mathcal U}\psi(x,Q)$, with $\mathcal U$ being a so-called uncertainty set or ambiguity set that is calibrated from data and, at least intuitively speaking, has a high likelihood of containing the true distribution.

Our main question to address is as follows: Considering the EO solution as a natural baseline, could we possibly improve its statistical performance in terms of the optimality gap \eqref{opt gap} (or equivalently the attained true objective value), by incorporating some regularization or robustness-enhancing modification?
Our main assertion is that, under standard conditions, it is theoretically \emph{impossible} to improve the statistical performance of EO in this regard. 

Our assertion is qualified by the following setup and conditions. First, we denote $\hat x_n^{EO}=x(\hat P_n)$ as the EO solution obtained from $n$ data points, where the subscript $n$ in $\hat x_n^{EO}$ and $\hat P_n$ highlights the dependence on sample size for the solution and the empirical distribution, and here $x(\cdot)$ is viewed as a function on $\hat P_n$. By a modification to EO, we mean to consider a wider choice of data-driven solution $\hat x_n^{\lambda}=x(\hat P_n,\lambda)$, where $\lambda$ is a tuning parameter in an expanded class of procedures that cover EO in particular. Without loss of generality, we set $x(\hat P_n,0)=x(\hat P_n)$, i.e., $\lambda=0$ corresponds to EO. This $\lambda$  appears virtually in all common regularization approaches, and also calibration proposals on the set size in DRO with a ``consistent" uncertainty set, i.e., a set with the property that it reduces to the singleton on the empirical distribution when its size $\lambda$ is tuned to be zero. Typically, $\lambda$ is chosen depending on the sample size $n$, and as $n$ grows, $\lambda$ eventually shrinks to $0$. For example, in regularization $\lambda$ represents a bias-variance tradeoff parameter introduced to avoid overfitting. When $n$ is large, the variance diminishes and so is the need to trade off a decrease in variance with an increase in bias, deeming a shrinkage of $\lambda$ to 0. Similarly, in DRO that uses a consistent confidence region as the uncertainty set, the set size converges to zero as $n$ gets large. Hereafter, when the dependence of $\lambda$ on $n$ is needed, we use the notation $\lambda_n$. 

Under the above setup, we have two main conditions. First, we focus on solutions $x(\cdot,\cdot)$ that are smooth with respect to the distribution $P$ and the tuning parameter $\lambda$. This condition is implied by the smoothness of the objective function $\psi(x,P)$ (with respect to both $x$ and $P$). Second, we consider the decision dimension $d$ to be fixed or, in other words, the setting where the sample size $n$ is large relative to the dimension. To summarize, we look at the most basic setting of smooth stochastic optimization in a large-sample regime, where our goal is to provide a fundamental argument to show the superiority of EO over all other possibilities.




We now explain the impossibility of statistical improvement. By this we mean that, in the large-sample regime, the optimality gap evaluated at the EO solution $\hat x_n^{EO}$ is always no worse than that evaluated at $\hat x_n^{\lambda_n}$, in terms of the risk profile measured by \emph{second-order stochastic dominance, regardless of how $\lambda_n$ depends on $n$}. 
As an immediate implication, this means for any non-decreasing convex function $f:\mathbb R\to\mathbb R$, we have
$$\mathbb E[f(\mathcal G(\hat x_n^{EO}))]\leq\mathbb E[f(\mathcal G(\hat x_n^{\lambda_n}))]$$
up to a negligible error, where $\mathbb E[\cdot]$ denotes the expectation with respect to the data used to obtain $\hat x_n^{EO}$ or $\hat x_n^{\lambda_n}$. In particular, setting $f(y)=y^2$, we conclude that the mean squared error of $Z(\hat x_n^{EO})$ against $Z(x^*)$ is always no larger than that of $Z(\hat x_n^{\lambda_n})$. The same conclusion holds for $f(y)=y^p$ for any other $p\geq1$. 


We will show how this impossibility result applies to all common examples of regularization on EO, and common DROs with a consistent uncertainty set. Moreover, we also show how our result applies to parametric settings, where the optimization involve unknown finite-dimensional parameters that need to be estimated. In the latter settings, the most straightforward approach is to plug in a consistent point estimate of the parameter into the optimization formulation, where this point estimate can be obtained by any common techniques such as maximum likelihood estimation (MLE) and the method of moments (MM) (\citealt{van2000asymptotic}, Chapters 4, 5). Here, we can again consider adding regularization on the plugged-in optimization formulation. We may also choose to use Bayesian approaches such as minimizing the expected posterior cost \citep{wu2018bayesian}, and consider regularizing the Bayesian formulation. In all the above cases, our result concludes that it is impossible to improve the solution obtained from EO or a simple consistent parameter estimation, in terms of the asymptotic risk profile of the optimality gap, by injecting regularization or DRO.


We will connect and contrast our results with other viewpoints in the data-driven optimization as well as the statistics literature. A reader who is proficient in stochastic optimization may find our claim on the superiority of EO very natural. Yet, as far as we know, our studied perspective appears unknown in the literature. Our results are intended to guide optimizers against engaging in suboptimal strategies in the considered basic large-sample situations, as we show in this situation that there is no theoretical improvement in using any strategies over EO. On the other hand, they are not intended to undermine the diverse strategies in data-driven optimization, as there are other situations where these strategies are used for good reasons (we will discuss these in Section \ref{sec:previous}).
 
In the following, Section \ref{sec:main} presents our main result, and Section \ref{sec:examples} discusses its applications on the range of data-driven formulations mentioned above. Then, in Section \ref{sec:previous}, we compare our results with established viewpoints in the DRO literature and classical statistics, and discuss scenarios beyond our considered setting in which alternate approaches to EO offer advantages. 

\section{Main Results}\label{sec:main}
Recall that $x^*$ is a minimizer of $Z(x)$. Also recall that $n$ is the sample size in an i.i.d. data set $\{\xi_i,i=1,\ldots,n\}$, and $\hat P_n$ denotes its empirical distribution, i.e., $\hat P(\cdot)=(1/n)\sum_{i=1}^n\delta_{\xi_i}(\cdot)$ where $\delta_{\xi_i}(\cdot)$ is the Dirac measure at the $i$-th data point $\xi_i$. $\hat x_n^{EO}=x(\hat P_n)$ is the EO solution, and $\hat x_n^{\lambda}=x(\hat P_n,\lambda)$ is the obtained solution from an expanded data-driven optimization procedure parametrized by $\lambda$ where $\hat x_n^0=\hat x_n^{EO}$. 

We set up some notation. In the following, we denote $E_Q[\cdot]$, $Var_Q(\cdot)$ and $Cov_Q(\cdot,\cdot)$ as the expectation, variance and covariance under a probability distribution $Q$. We use ``$\Rightarrow$" to denote weak convergence or convergence in distribution, ``$\stackrel{p}{\to}$" denote convergence in probability,  ``$\stackrel{d}{=}$" denote equality in distribution, and ``a.s." denote almost surely. We denote $\|\cdot\|$ as the Euclidean norm. For any deterministic sequences $a_k\in\mathbb R$ and $b_k\in\mathbb R$, both indexed by a common index, say $k$ that goes to $\infty$, we say that $a_k=o(b_k)$ if $a_k/b_k\to0$, $a_k=O(b_k)$ if there exists a finite $M>0$ such that $|a_k/b_k|<M$ for all sufficiently large $k$, $a_k=\omega(b_k)$ if $|a_k/b_k|\to\infty$, $a_k=\Omega(b_k)$ if there exists a finite $M>0$ such that $|a_k/b_k|>M$ for all sufficiently large $k$, and $a_k=\Theta(b_k)$ if there exist finite $\underline M, \overline M>0$ such that $\underline M<|a_k/b_k|<\overline M$ for all sufficiently large $k$. We use the notations $o_p(\cdot)$ and $O_p(\cdot)$ to denote a smaller and an at most equal stochastic order respectively. Namely, for a sequence of  random vector $A_k\in\mathbb R^d$ and deterministic sequence $b_k\in\mathbb R$, both indexed by a common index $k$ that goes to $\infty$, $A_k=o_p(b_k)$ means $A_k/b_k\stackrel{p}{\to}0$ as $k\to\infty$. Correspondingly, $A_k=O_p(b_k)$ means for any $\epsilon>0$, there exists a large enough $N>0$ and $M>0$ such that $P(\|A_k/b_k\|\leq M)\geq1-\epsilon$ for any $k>N$. Finally, for a sequence of random variable $A_k\in\mathbb R$, we say that $A_k\stackrel{p}{\to}\infty$ if for any $\epsilon>0$ and $M>0$ there exists $N>0$ large enough such that $P(A_k>M)\geq1-\epsilon$ for any $k>N$. Similarly, $A_k\stackrel{p}{\to}-\infty$ if for any $\epsilon>0$ and $M>0$ there exists $N>0$ large enough such that $P(A_k<-M)\geq1-\epsilon$ for any $k>N$.  Finally, we denote $\text{ess sup}_Q\ f$ as the essential supremum of a random function $f$ under distribution $Q$, and $\top$ as the transpose.

\subsection{Conditions}
We make three assumptions:
\begin{assumption}[Optimality conditions]
A true minimizer $x^*$ for $Z(x)$ satisfies the second-order optimality conditions, namely the gradient $\nabla Z(x^*)=0$ and the Hessian $\nabla^2Z(x^*)$ is positive semidefinite.\label{optimality}
\end{assumption}

\begin{assumption}[Linearizability]
The data-driven solution $\hat x_n^\lambda$ satisfies
\begin{equation}\hat x_n^\lambda-x^*=\langle IF(\xi),\hat P_n-P\rangle+\lambda K+o_p\left(\frac{1}{\sqrt n}+\lambda\right)\label{asymptotic expansion}
\end{equation}
as $n\to\infty$ and $\lambda\to0$ (at any rate), for some function $IF:\Xi\to\mathbb R^d$ and constant vector $K\in\mathbb R^d$. Moreover, $Cov_P(IF(\xi))$, the covariance matrix of $IF(\xi)$ under $P$, is entry-wise finite. Furthermore, when $\lambda=0$, we have $\hat x_n^0=\hat x_n^{EO}$. \label{solution expansion}
\end{assumption}

\begin{assumption}[Non-degeneracy]
We have $K^\top\nabla^2Z(x^*)K\neq0$ where $\nabla^2Z(x^*)$ and $K$ are defined in Assumptions \ref{optimality} and \ref{solution expansion} respectively.\label{nondegenerate}
\end{assumption}

Assumption \ref{optimality} is the standard optimality condition for unconstrained problem or when $x^*$ is in the interior of the feasible region. Note that we implicitly assume that $Z(\cdot)$ is twice differentiable in the assumption. The unconstrained or interior-point nature of the solution is not necessary and is used here to present a cleaner setting (we will extend our result later to the constrained case). 

Assumption \ref{solution expansion} is a smoothness condition on $\hat x_n^\lambda=x(\hat P_n,\lambda)$ with respect to both $\hat P_n$ and $\lambda$. In particular, it states that $\hat x_n^\lambda$ converges to $x^*$ linearly in $\hat P_n-P$ as $n\to\infty$ and $\lambda$ as $\lambda\to0$, with a negligible remainder term of higher order than $1/\sqrt n+\lambda$. More specifically, the function $IF(\cdot)$ in \eqref{asymptotic expansion} is known as the influence function \citep{hampel1974influence} at $P$, which is viewed as the functional gradient of $x(\cdot,\lambda)$ with respect to the probability distribution $P$. This gradient is a function on the image of the random object generated by $P$, i.e., the dual space of $P$, and $\langle\cdot,\cdot\rangle$ is an inner product that is written more explicitly as
\begin{equation}
\langle IF(\xi),\hat P_n-P\rangle:=\int IF(\xi)d(\hat P_n-P)=\int IF(\xi)d\hat P_n=\frac{1}{n}\sum_{i=1}^nIF(\xi_i)\label{IF}
\end{equation}
Here in \eqref{IF}, the first equality is the definition of the inner product $\langle\cdot,\cdot\rangle$. The second equality follows by noting that we can always assume $IF(\xi)$ satisfies $E_P[IF(\xi)]=0$ (because $E_{\hat P_n}[1]=1$, we can take $IF(\xi)-E_P[IF(\xi)]$ to be our new influence function if not the case). The third equality follows from the definition of the empirical distribution $\hat P_n$. Note that when $\lambda=0$, Assumption \ref{solution expansion} implies that the EO solution $\hat x_n^{EO}=\hat x_n^0$ satisfies
\begin{equation}
\hat x_n^{EO}-x^*=\langle IF(\xi),\hat P_n-P\rangle+o_p\left(\frac{1}{\sqrt n}\right)\label{asymptotic expansion original}
\end{equation}
Such a linear relation can in fact arise more generally, for instance under Hadamard differentiability (\citealt{van2000asymptotic} Chapter 20), though this is not needed for our current purpose.


Note that Assumption \ref{solution expansion} also stipulates that $K$ is the (partial) derivative of $x(P,\lambda)$ with respect to $\lambda$. Moreover, \eqref{asymptotic expansion} and \eqref{asymptotic expansion original} imply the consistency of solutions $\hat x_n^\lambda$ and $\hat x_n$, in the sense that $\hat x_n^\lambda\stackrel{p}{\to}x^*$ and $\hat x_n^{EO}\stackrel{p}{\to}x^*$ as $n\to\infty$ and $\lambda\to0$ (at any rate). 

Finally, Assumption \ref{nondegenerate} is a technical non-degeneracy assumption to guarantee that the impact due to the introduction of $\lambda$ in $\hat x_n^\lambda$ is not completely canceled out by the Hessian $\nabla Z^2(x^*)$. Note that by Assumption \ref{optimality} we always have $K^\top\nabla^2Z(x^*)K\geq0$, and Assumption \ref{nondegenerate} guarantees that this inequality is strict. 



\subsection{Second-Order Stochastic Dominance}
To explain our results, we need the following concept \citep{rothschild1970increasing}:
\begin{definition}[Second-order stochastic dominance for losses]
For any two real-valued random variables $A$ and $B$, we say that $A$ is second-order stochastically dominated by $B$ if $E[u(A)]\leq E[u(B)]$ for any non-decreasing convex function $u:\mathbb R\to\mathbb R$.\label{sosd}
\end{definition}

Definition \ref{sosd} is also known as a monotone convex order, in which $A$ is said to be smaller than $B$ in increasing convex order (\citealt{shaked2007stochastic} Section 4.A.1). The following is an important equivalence relation:
\begin{proposition}[Decomposition via mean-preserving spread]
For any two real-valued random variables $A$ and $B$, we have that $A$ is second-order stochastically dominated by $B$ if and only if
\begin{equation}
\tilde B=\tilde A+\eta+\epsilon\label{decomposition}
\end{equation}
for some $\tilde A\stackrel{d}{=}A$, $\tilde B\stackrel{d}{=}B$, an a.s. non-negative random variable $\eta$, and a random variable $\epsilon$ such that $\tilde B$ is a mean-preserving spread of $\tilde A+\eta$, i.e.,
$$E[\epsilon|\tilde A+\eta]=0\text{\ \  a.s.}$$\label{mean-preserving spread}
\end{proposition}

Note that the definition of second-order stochastic dominance is defined solely on the probability distributions of $A$ and $B$ above. It is common to define this notion on the distributions directly, though here we use random variable for the convenience of our subsequent developments. Moreover, let us make clear that our Definition \ref{sosd} applies to loss (i.e., smaller is desirable) rather than gain (i.e., bigger is desirable), the latter more customarily used in the economics literature \citep{hanoch1969efficiency,hadar1969rules,rothschild1970increasing}. This distinction, however, is immaterial, as we can simply view gain as negative loss. That is, if we call the gain $-A$ to second-order stochastically dominate the gain $-B$, then Definition \ref{sosd} states that $E[-u(-(-A))]\geq E[-u(-(-B))]$ for any non-decreasing convex function $u$, which is equivalent to $E[v(-A)]\geq E[v(-B)]$ for any non-decreasing concave function $v$. This reduces back to notion that a gain is preferable if it has a higher expected utility.

Proposition \ref{mean-preserving spread} is well-established (e.g., \citealt{shaked2007stochastic} Theorem 4.A.5). Proposition \ref{mean-preserving spread} states that a second-order stochastically dominant variable $B$, when compared to $A$, contains two additional terms $\eta$ and $\epsilon$. The first term is a non-negative random variable $\eta$, so that $A+\eta$ is less attractive than $A$ in terms of first-order stochastic dominance, i.e., the distribution function of $A+\eta$ is at least that of $A$. The second term is $\epsilon$ that does not change the expectation but adds more variability to $A+\eta$. This means that $B$ is a so-called \emph{mean-preserving spread} of $A+\eta$ \citep{landsberger1993mean}. From a risk-averse perspective, additional uncertainty caused by higher variability is always undesirable. Second-order stochastic dominance thus stipulates that $B$ is less desirable than $A$ in terms of both the notions that ``smaller is desirable" and ``less uncertainty is desirable".

Finally, we introduce the following notion that facilitates the use of stochastic dominance for weak limits, which occurs in our studied large-sample regime:
\begin{definition}[Asymptotic second-order stochastic dominance]
For any two real-valued random sequences $A_n$ and $B_n$, we say that $A_n$ is asymptotically second-order stochastically dominated by $B_n$ as $n\to\infty$ if $A_n\Rightarrow W_0$ and $B_n\Rightarrow W_1$ such that $W_0$ is second-order stochastically dominated by $W_1$ or, more generally, for any subsequences $n_{k_0},n_{k_1}\to\infty$ such that $A_{n_{k_0}}$ and $B_{n_{k_1}}$ have weak limits, say $W_0$ and $W_1$ respectively, $W_0$ is second-order stochastically dominated by $W_1$.\label{asymptotic sosd}
\end{definition}

``Asymptotic" here means that the second-order stochastic dominance holds in the weak limits or, in the case that no unique weak limits exist, then we look at all possible subsequences that have weak limits. We also remark that the weak convergence here is broadly defined as including the case where the limits $W_0$ and $W_1$ can be $\pm\infty$, in which case it means the corresponding sequences $\stackrel{p}{\to}\pm\infty$.

\subsection{The Impossibility Theorem}
With the above setting and notation, we are now ready to state our main result. 

\begin{theorem}[Impossibility of statistical improvement from EO]
Suppose Assumptions \ref{optimality}, \ref{solution expansion} and \ref{nondegenerate} hold, and $\lambda_n\to0$ as $n\to\infty$. Then the optimality gap of $\hat x_n^{EO}$ scaled by $n$, namely $n\mathcal G(\hat x_n^{EO})$, is asymptotically second-order stochastically dominated by that of $\hat x_n^{\lambda_n}$, namely $n\mathcal G(\hat x_n^{\lambda_n})$.\label{main thm}
\end{theorem}


Theorem \ref{main thm} states that no matter how we choose $\lambda_n$ in the expanded data-driven procedure, as long as it goes to 0, then it is impossible to improve $\hat x_n^{EO}$ asymptotically in terms of the risk profile of the optimality gap measured by second-order stochastic dominance. Note that the scaling $n$ in front of the optimality gap in the theorem is natural as this is the scaling that gives a nontrivial limit under the first-order optimality condition.

The key in showing Theorem \ref{main thm} is to decompose the weak limit of $n\mathcal G(\hat x_n^{\lambda_n})$ into that of $n\mathcal G(\hat x_n^{EO})$, a positive term, and a conditionally unbiased noise term that leads to the mean-preserving spread as in \eqref{decomposition}. Essentially, the two extra terms signify that introducing $\lambda_n$ would simultaneously lead to an extra ``bias" and an extra ``variability". This behavior holds regardless of how we choose the sequence $\lambda_n$. The detailed proof of Theorem \ref{main thm} is in Appendix \ref{sec:proof main}, and the next subsection gives the roadmap and intuitive explanation.

\subsection{Explaining the Impossibility Theorem}\label{sec:development}
To explain Theorem \ref{main thm} in more detail, we start with writing the mathematical form of the optimality gap more explicitly as follows:
\begin{lemma}
Under Assumptions \ref{optimality} and \ref{solution expansion}, as $n\to\infty$ and $\lambda\to0$, the optimality gap of $\hat x_n^\lambda$ satisfies
\begin{eqnarray}
\mathcal G(\hat x_n^\lambda)&:=&Z(\hat x_n^\lambda)-Z(x^*)\notag\\
&=&\frac{1}{2}\langle IF(\xi),\hat P_n-P\rangle^\top\nabla^2Z(x^*)\langle IF(\xi),\hat P_n-P\rangle+\frac{1}{2}\lambda^2K^\top\nabla^2Z(x^*)K+\lambda K^\top\nabla^2Z(x^*)\langle IF(\xi),\hat P_n-P\rangle{}\notag\\
&&{}+o_p\left(\frac{1}{n}+\lambda^2\right)\label{prelim expression}
\end{eqnarray}\label{explicit lemma}
\end{lemma}

Lemma \ref{explicit lemma} can be obtained straightforwardly by Taylor-expanding $Z(\hat x_n^\lambda)$ around $x^*$, giving
\begin{equation}
Z(\hat x_n^\lambda)-Z(x^*)=\frac{1}{2}(\hat x_n^\lambda-x^*)^\top\nabla^2Z(x^*)(\hat x_n^\lambda-x^*)+o_p(\|\hat x_n^\lambda-x^*\|)\label{new opt gap CR}
\end{equation}
and plugging in the linear approximation of $\hat x_n^\lambda-x^*$ in \eqref{asymptotic expansion}. In \eqref{prelim expression}, the first term is approximately $\mathcal G(\hat x_n^{EO})$ (i.e., when setting $\lambda=0$), while the second and third term appear as the additional quadratic and cross terms arising from the ``bias" $\lambda K$ in \eqref{asymptotic expansion} for $\hat x_n^\lambda$. In other words, Lemma \ref{explicit lemma} deduces the following relation between the optimality gaps of $\hat x_n^\lambda$  and $\hat x_n^{EO}$
$$\mathcal G(\hat x_n^\lambda)\stackrel{d}{=}\mathcal G(\hat x_n^{EO})+\frac{1}{2}\lambda^2K^\top\nabla^2Z(x^*)K+\lambda K^\top\nabla^2Z(x^*)\langle IF(\xi),\hat P_n-P\rangle{}+o_p\left(\frac{1}{n}+\lambda^2\right)$$
or equivalently the relation between the true objective values attained by $\hat x_n^\lambda$  and $\hat x_n^{EO}$
\begin{equation}
Z(\hat x_n^\lambda)\stackrel{d}{=}Z(\hat x_n^{EO})+\frac{1}{2}\lambda^2K^\top\nabla^2Z(x^*)K+\lambda K^\top\nabla^2Z(x^*)\langle IF(\xi),\hat P_n-P\rangle{}+o_p\left(\frac{1}{n}+\lambda^2\right)\label{obj value comparison}
\end{equation}
The detailed proof of Lemma \ref{explicit lemma} is left to Appendix \ref{sec:proof main}. 

Now we highlight the main intuition in obtaining Theorem \ref{main thm}. By the definition of $\langle IF(\xi),\hat P_n-P\rangle$ in \eqref{IF}, the central limit theorem (CLT) implies that
$$\langle IF(\xi),\hat P_n-P\rangle\approx\frac{Y}{\sqrt n}$$
in distribution, where $Y$ is a Gaussian vector with mean 0 and covariance matrix $Cov_P(IF(\xi))$. Thus we can write the expression of $\mathcal G(\hat x(\lambda))$ in Lemma \ref{explicit lemma} as
\begin{equation}
\mathcal G(\hat x_n^\lambda)
\approx\frac{1}{2n}Y^\top\nabla^2Z(x^*)Y+\frac{1}{2}\lambda^2K^\top\nabla^2Z(x^*)K+\frac{1}{\sqrt n}\lambda K^\top\nabla^2Z(x^*)Y+o_p\left(\frac{1}{n}+\lambda^2\right)\label{interim6}
\end{equation}
in distribution. Moreover, setting $\lambda=0$, the EO optimality gap becomes
\begin{equation}
\mathcal G(\hat x_n^{EO})
\approx\frac{1}{2n}Y^\top\nabla^2Z(x^*)Y+o_p\left(\frac{1}{n}\right)\label{interim7}
\end{equation}

Now consider all possible choices of the sequence $\lambda$ in relation to $n$, and for each case we see how the first three terms in \eqref{interim6} behaves. Suppose $\lambda_n=o(1/\sqrt n)$. Then the second and third terms are of smaller order than the first term, so that \eqref{interim6} reduces to $(1/(2n))Y^\top\nabla^2Z(x^*)Y+o_p(1/n)$. In this case, $\mathcal G(\hat x_n^{\lambda_n})$ and $\mathcal G(\hat x_n^{EO})$ behave the same asymptotically. In other words, the expanded procedure does not offer any first-order benefit relative to EO. Now suppose, on the other hand, that $\lambda_n=\omega(1/\sqrt n)$. Then the second and third terms are both of bigger order than the first term, with the second term the most dominant, and so \eqref{interim6} becomes $\frac{1}{2}\lambda_n^2K^\top\nabla^2Z(x^*)K+o_p\left(\lambda_n^2\right)$. In this case, $\mathcal G(\hat x_n^{\lambda_n})$ is of bigger order than $\mathcal G(\hat x_n^{EO})$, so that the expanded procedure gives a solution that is worse than EO. Thus, we are left with choosing $\lambda_n=\Theta(1/\sqrt n)$. 

Suppose $\lambda_n\approx a/\sqrt n$ for some finite $a\neq0$ as $n\to\infty$. We have
$$\mathcal G(\hat x_n^{\lambda_n})
\approx\frac{1}{2n}Y^\top\nabla^2Z(x^*)Y+\frac{1}{2n}a^2K^\top\nabla^2Z(x^*)K+\frac{1}{n}a K^\top\nabla^2Z(x^*)Y+o_p\left(\frac{1}{n}\right)$$
so that all three terms have the same order $1/n$. The coefficient in this first-order term is 
\begin{equation}
\frac{1}{2}Y^\top\nabla^2Z(x^*)Y+\frac{1}{2}a^2K^\top\nabla^2Z(x^*)K+a K^\top\nabla^2Z(x^*)Y\label{interim8}
\end{equation}
Note that the second term in \eqref{interim8} is deterministic and always non-negative. On the other hand, the third term has conditional mean zero given the first term, namely
\begin{equation}
E\left[a K^\top\nabla^2Z(x^*)Y\Bigg|\frac{1}{2}Y^\top\nabla^2Z(x^*)Y\right]=0\label{interim9}
\end{equation}
To see this, note that since $\nabla^2Z(x^*)$ is positive semidefinite by Assumption \ref{optimality}, $Y^\top\nabla^2Z(x^*)Y$ can be written as a sum-of-squares of linear transformations of $Y$, i.e., $Y^\top\nabla^2Z(x^*)Y=\|(\nabla^2Z(x^*))^{1/2}Y\|^2$ where $(\nabla^2Z(x^*))^{1/2}$ denotes the square-root matrix of $\nabla^2Z(x^*)$. Moreover, note that $Y$, as a mean-zero Gaussian vector, is symmetric (i.e., the densities at $y$ and $-y$ are the same). Thus, conditional on the knowledge of $\frac{1}{2}Y^\top\nabla^2Z(x^*)Y$, it is equally likely for $Y_j$ to take $y_j$ and $-y_j$ for any $y_j$, and we have $E[Y_j|Y^\top\nabla^2Z(x^*)Y]=0$
for all components $Y_j,j=1,\ldots,d$ of $Y$. This implies \eqref{interim9}. In other words, \eqref{interim8} is a mean-preserving spread of $(1/2)Y^\top\nabla^2Z(x^*)Y+(1/2)a^2K^\top\nabla^2Z(x^*)K$.

Putting the above together, we see that, in the case $\lambda_n\approx a/\sqrt n$ for finite $a\neq0$, the first-order term of $\mathcal G(\hat x_n^{\lambda_n})$, \eqref{interim8}, is exactly decomposable into the form in \eqref{decomposition} in Proposition \ref{mean-preserving spread}. Compared with the first-order coefficient of $\mathcal G(\hat x_n^{EO})$, namely $(1/2)Y^\top\nabla^2Z(x^*)Y$, we thus have $\mathcal G(\hat x_n^{EO})$ second-order stochastically dominated by $\mathcal G(\hat x_n^{\lambda_n})$ in terms of first-order behavior. Therefore, in all possible choices of $\lambda_n$, the asymptotic statistical behavior of $\mathcal G(\hat x_n^{\lambda_n})$ cannot be better than $\mathcal G(\hat x_n^{EO})$. 


We summarize the above intuitive explanation with the following two propositions:
\begin{proposition}[A trichotomy of optimality gap behaviors]
Under Assumptions \ref{optimality}, \ref{solution expansion} and \ref{nondegenerate}, suppose $\lambda_n\to0$ with $\sqrt n\lambda_n\to a$ for some $a\in\mathbb R\cup\{-\infty,\infty\}$ as $n\to\infty$. We have the following regarding the optimality gap of $\hat x_n^{\lambda_n}$:
\begin{enumerate}
\item When $a=0$, then 
$$n\mathcal G(\hat x_n^{\lambda_n})\Rightarrow\frac{1}{2}Y^\top\nabla^2Z(x^*)Y$$
\item When $a=\infty$ or $-\infty$, then 
$$\frac{1}{\lambda_n^2}\mathcal G(\hat x_n^{\lambda_n})\Rightarrow\frac{1}{2}K^\top\nabla^2Z(x^*)K$$
\item When $0<|a|<\infty$, then
$$n\mathcal G(\hat x_n^{\lambda_n})\Rightarrow\frac{1}{2}Y^\top\nabla^2Z(x^*)Y+\frac{1}{2}a^2K^\top\nabla^2Z(x^*)K+aK^\top\nabla^2Z(x^*)Y$$
\end{enumerate}
On the other hand, the optimality gap of $\hat x_n$ satisfies, as $n\to\infty$,
$$n\mathcal G(\hat x_n^{EO})\Rightarrow\frac{1}{2}Y^\top\nabla^2Z(x^*)Y$$
In all the above, $Y$ denotes a Gaussian random vector with mean 0 and covariance matrix $Cov_P(IF(\xi))$.
\label{trichotomy}
\end{proposition}


\begin{proposition}[Comparisons on the trichotomy]
Under the same assumptions and notations in Proposition \ref{trichotomy}, consider $n\to\infty$. When $a=0$, $n\mathcal G(\hat x_n^{\lambda_n})$ and $n\mathcal G(\hat x_n^{EO})$ have the same weak limit. When $a=\infty$ or $-\infty$, $n\mathcal G(\hat x_n^{\lambda_n})\stackrel{p}{\to}\infty$ whereas $n\mathcal G(\hat x_n^{EO})$ converges weakly to a tight random variable. When $0<|a|<\infty$, the weak limit of $n\mathcal G(\hat x_n^{EO})$ is second-order stochastically dominated by that of $n\mathcal G(\hat x_n^{\lambda_n})$.\label{comparisons}
\end{proposition}

Proposition \ref{trichotomy} summarizes the asymptotic limits of $\mathcal G(\hat x_n^{\lambda_n})$ derived from Lemma \ref{explicit lemma}. Proposition \ref{comparisons} then compares them with that of $\mathcal G(\hat x_n^{EO})$, and concludes that asymptotic second-order stochastic dominance holds in all cases, thus concluding Theorem \ref{main thm}.

We close this section by cautioning that when Assumption \ref{nondegenerate} is violated, so that $K^\top\nabla^2Z(x^*)K=0$, then we would have $(\nabla^2Z(x^*))^{1/2}K=0$ and hence $K^\top\nabla^2Z(x^*)=0$. Thus both the second and third terms in \eqref{interim8} vanish. In this case, choosing a suitable $\lambda_n=\omega(1/\sqrt n)$ potentially retains the first-order behavior of $\mathcal G(\hat x_n^{\lambda_n})$ to be comparable to $\mathcal G(\hat x_n^{EO})$, i.e., of order $1/n$, with an involved coefficient that could be bigger or smaller. Thus in this degenerate case the comparison is inconclusive.

\subsection{Generalization to Constrained Problems}\label{sec:constrained}
We close this section by presenting the generalization of our results to constrained optimization. Consider
\begin{equation}
\min_x\ Z(x)\ \ \text{subject to}\ \ g_j(x)\leq0\text{\ for\ }j\in J_1,\ \ g_j(x)=0\text{\ for\ }j\in J_2\label{opt con}
\end{equation}
That is, we have $|J_1|$ inequality constraints and $|J_2|$ equality constraints, with constraint functions denoted $g_j(x)$. We could also add non-negativity constraints on $x$, as long as the solution satisfies the assumptions we make momentarily. We assume the following first-order optimality conditions:
\begin{assumption}[Lagrangian optimality condition]
The optimal solution $x^*$ to \eqref{opt con} solves
$$\nabla Z(x)+\sum_{j\in B}\alpha_j^*\nabla g_j(x)=0$$
where $\alpha_j^*$'s are the Lagrange multipliers, and $B$ indicates the binding set of constraints, i.e., $B=\{j\in J_1\cup J_2:g_j(x^*)=0\}$. Moreover, $\nabla Z^2(x)+\sum_{j\in B}\alpha_j^*\nabla^2g_j(x)=0$ is positive semidefinite.
\label{assumption:Lagrangian}
\end{assumption}

The first condition in Assumption \ref{assumption:Lagrangian} is the standard KKT condition. The second condition is the second-order optimality condition imposed on the Lagrangian. Note that we have implicitly assumed $Z$ and $g_j$'s are twice differentiable in the assumption.

In addition, we assume the following behavior on $\hat x_n^\lambda$:
\begin{assumption}[Preservance of binding constraints]
The data-driven solution $\hat x_n^\lambda$ satisfies $g_j(\hat x_n^\lambda)=0$ for any $j\in B$, with probability converging to 1 as $n\to\infty$ and $\lambda\to0$.\label{assumption:preservance}
\end{assumption}

Assumption \ref{assumption:preservance} means that, asymptotically, the data-driven solution retains the same set of binding constraints as the true solution. This condition typically holds for any solution that is consistent in converging to $x^*$. Next, in parallel to Assumption \ref{nondegenerate}, we have the following non-degeneracy condition in the constrained case:
\begin{assumption}[Non-degeneracy in constrained problems]
We have $K^\top(\nabla^2Z(x^*)+\sum_{j\in B}\alpha_j^*\nabla^2g_j(x^*))K\neq0$ where $\nabla^2Z(x^*)$, $\alpha_j^*$, $\nabla^2g_j(x^*)$ and $K$ are defined in Assumptions \ref{assumption:Lagrangian} and \ref{solution expansion}.\label{assumption:nondegenerate con}
\end{assumption}

With the above assumptions, we now argue that all the impossibility results we have discussed hold in the constrained setting, as follows:
\begin{theorem}[Impossibility of statistical improvement over EO in constrained problems]
Suppose Assumptions \ref{assumption:Lagrangian}, \ref{assumption:preservance} and \ref{assumption:nondegenerate con} hold, and so does the expansion of solution $\hat x_n^\lambda$ depicted in Assumption \ref{solution expansion}. Then the conclusions of Theorem \ref{main thm}, Lemma \ref{explicit lemma}, Proposition \ref{trichotomy} and Proposition \ref{comparisons} all hold, except that $\nabla^2Z(x^*)$ is replaced by $\nabla^2Z(x^*)+\sum_{j\in B}\alpha_j^*\nabla^2g_j(x^*)$ in all these results.\label{main thm con}
\end{theorem}

The proof of Theorem \ref{main thm con} follows a similar development as in Section \ref{sec:development}, but instead of using the unconstrained first-order optimality to remove the first-order term in the Taylor series expansion of $\mathcal G(\hat x_n^\lambda)$, we use the Lagrangian (Assumption \ref{assumption:Lagrangian}) to express this term in terms of the constraints, which are in turn analyzed via another Taylor expansion. The detailed proof is in Appendix \ref{sec:proof main}.

\section{Applications to Data-Driven Optimization Formulations}\label{sec:examples}
We apply Theorem \ref{main thm} to several common data-driven optimization formulations, including regularization (Section \ref{sec:regularization}), DRO (Section \ref{sec:DRO}), parametric optimization (Section \ref{sec:parametric}) and its Bayesian counterpart (Section \ref{sec:Bayesian}). For each data-driven optimization approach, we will identify the expansion \eqref{asymptotic expansion} in Assumption \ref{solution expansion}, the main condition to achieve our impossibility result. This requires certain regularity conditions which for succinctness are presented in the Appendix.

\subsection{Regularization}\label{sec:regularization}
Consider $Z(x):=\psi(x,P)$ as the expected value objective function $E_P[h(x,\xi)]$ where $h:\mathcal X\times\Xi\to\mathbb R$. The EO objective function is then $\psi(x,\hat P_n)=E_{\hat P_n}[h(x,\xi)]$. We consider $\hat x_n^{Reg,\lambda}$ as a regularized solution obtained from 
\begin{equation}
\min_{x\in\mathcal X}\{\psi(x,\hat P_n)+\lambda R(x)=E_{\hat P_n}[h(x,\xi)]+\lambda R(x)\}\label{reg def}
\end{equation}
where $R:\mathcal X\to\mathbb R$ is a penalty function on $x$ and $\lambda\geq0$ is the regularization parameter. For instance, when $h$ represents the least-square objective or some loss function and $R(x)=\|x\|^2$, we have a ridge regression (e.g., \citealt{friedman2001elements}).


We have the following result:
\begin{theorem}[Impossibility of improving EO via regularization]
Consider $Z(x)=E_P[h(x,\xi)]$ with a minimizer $x^*$. Suppose we collect i.i.d. observations $\xi_1,\ldots,\xi_n$ with empirical distribution $\hat P_n$, and let $\hat x_n^{EO}$ be a minimizer of $E_{\hat P_n}[h(x,\xi)]$ and $\hat x_n^{Reg,\lambda}$ a minimizer of \eqref{reg def}. Under Assumptions \ref{regularity loss}--\ref{first-order assumption} in Appendix \ref{sec:proof regularization}, we have
\begin{eqnarray}
&&\hat x_n^{Reg,\lambda_n}-x^*\notag\\
&=&-\langle (E_P[\nabla^2h(x^*,\xi)])^{-1}\nabla h(x^*,\xi),\hat P_n-P\rangle-\lambda_n(E_P[\nabla^2h(x^*,\xi)])^{-1}\nabla R(x^*)+o_p\left(\frac{1}{\sqrt n}+\lambda_n\|\nabla R(x^*)\|\right)\label{regularization expression}
\end{eqnarray}
as $n\to\infty$ and $\lambda_n\to0$ (at any rate relative to $n$). Hence, if in addition Assumptions \ref{optimality} and \ref{nondegenerate} hold with $K=-(E_P[\nabla^2h(x^*,\xi)])^{-1}\nabla R(x^*)$, then $n\mathcal G(\hat x_n^{EO})$ is asymptotically second-order stochastically dominated by $n\mathcal G(\hat x_n^{Reg,\lambda_n})$. 
\label{regularization thm}
\end{theorem}

The proof of Theorem \ref{regularization thm} is in Appendix \ref{sec:proof regularization}, which requires the theory of $M$-estimation with nusiance parameter.  In addition to smoothness  (Assumptions \ref{regularity loss}, \ref{regularity penalty}) and first-order optimality conditions (Assumption \ref{first-order assumption}), we also need to control the function complexity of the loss function \citep{van1996weak} via Donsker and Glivenko-Cantelli properties (Assumption \ref{Donsker}). Theorem \ref{regularization thm} can also be translated to DRO based on the Wasserstein ball thanks to its connection with regularization; see the next subsection.

At a first glance, Theorem \ref{regularization thm} may look contradictory to the regularization literature. For instance, even in the linear regression, it is known that a ridge regression with a properly chosen $\lambda$ can always improve the estimation quality under large sample \citep{li1986asymptotic,li1987asymptotic}. The catch here is the criterion to measure the quality of solution. Our considered criterion in this paper is the risk profile of the entire distribution of the optimality gap, whereas the criterion looked at in the ridge regression can be seen to correspond to the \emph{expected} optimality gap. This latter criterion appears reasonable for estimation problems, but it is not sufficient from an optimization viewpoint: An improvement in the expected gap does not necessarily mean the obtained solution is better statistically in terms of the attained objective value, as the variability of the optimality gap can become worse -- In fact this fundamental dilemma is precisely what we showed.

To distinguish clearly optimality gap versus expected optimality gap, consider the approximation of $\mathcal G(\hat x_n^\lambda)$ in \eqref{prelim expression}. Under standard regularity conditions, its expected value, taken with respect to all the data, is
\begin{eqnarray}
\mathbb E[\mathcal G(\hat x_n^\lambda)]&=&\mathbb E[Z(\hat x_n^\lambda)-Z(x^*)]\notag\\
&=&\mathbb E\Bigg[\frac{1}{2}\langle IF(\xi),\hat P_n-P\rangle^\top\nabla^2Z(x^*)\langle IF(\xi),\hat P_n-P\rangle+\frac{1}{2}\lambda^2K^\top\nabla^2Z(x^*)K{}\notag\\
&&{}+\lambda K^\top\nabla^2Z(x^*)\langle IF(\xi),\hat P_n-P\rangle\Bigg]+o\left(\frac{1}{n}+\lambda^2\right)\notag\\
&=&\mathbb E[\mathcal G(\hat x_n^{EO})]+\frac{1}{2}\lambda^2K^\top\nabla^2Z(x^*)K+o\left(\frac{1}{n}+\lambda^2\right)\label{expected gap}
\end{eqnarray}
where the last equality follows by comparing with the approximation for $\mathbb E[\mathcal G(\hat x_n^{EO})]$, and noting that $\frac{1}{2}\lambda^2K^\top\nabla^2Z(x^*)K$ is deterministic and $\lambda K^\top\nabla^2Z(x^*)\langle IF(\xi),\hat P_n-P\rangle$ has mean 0. The dominant term in the remainder, namely $o(1/n+\lambda^2)$, is actually of order $O(\lambda/n)$, which comes from the cross-term of two ``$\hat P_n-P$" and one ``$\lambda$" (the other higher-order terms either have expectation 0 or are dominated by others). Thus, we can choose $\lambda=k/n$, for some constant $k$, such that the second term in \eqref{expected gap} and this $O(\lambda/n)$ term possess the same $1/n^2$ order and are in opposite signs, thus leading to a smaller expected optimality gap over $\hat x_n^{EO}$. However, if the variability of the optimality gap is taken into account, then we need to consider the third term of \eqref{prelim expression} which is of stochastic order $\lambda/\sqrt n$. Choosing $\lambda=k/n$ then entails this term to order $1/n^{3/2}$, which is larger than the improvement of order $1/n^2$ in the expected gap and as a result washes away this gain.

We should also make clear that our impossibility result in Theorem \ref{main thm} does apply to the comparison in terms of the expected optimality gap, because the second-order stochastic dominance implies as a particular case that the expected values follow the corresponding ordering. This may appear again as a contradiction to our result in saying that regularization can improve the expected gap of EO. However, this latter gain is actually of higher-order than $1/n$, whereas in the dominant term of order $1/n$ in the optimality gap there is no gain, which is what our result implies. This means that, at least in low-dimensional cases, the improvement using regularization, even focusing only on the expected optimality gap, is negligible.

Finally, we justify one of our claims above that the criterion used in justifying the gain in ridge regression corresponds to the expected optimality gap. Consider the simple least-square problem $\min_{\beta\in\mathbb R^d} E(Y-X^\top\beta)^2$, where $\xi=(X,Y)\in\mathbb R^{d+1}$ follows a linear model $Y=X^\top\beta+\epsilon$ with $E[\epsilon|X]=0$, and $E[XX^\top]=I$. Suppose from data we obtain a coefficient estimate $\hat\beta$. Then the optimality gap $\mathcal G(\hat\beta)$, according to the least-square objective function, is
\begin{eqnarray*}
&&E(Y-X^\top\hat\beta)^2-E(Y-X^\top\beta)^2\\
&=&-2E[YX^\top(\hat\beta-\beta)]+E[\hat\beta^\top XX^\top\hat\beta-\beta^\top XX^\top\beta]\\
&=&-2E[YX^\top](\hat\beta-\beta)+(\hat\beta^\top\hat\beta-\beta^\top\beta)\\
&=&-2\beta^\top(\hat\beta-\beta)+\|\hat\beta\|^2-\|\beta\|^2\\
&=&\|\hat\beta-\beta\|^2
\end{eqnarray*}
where we have used $E[YX^\top]=\beta^\top$. That is, $\|\hat\beta-\beta\|^2$ is the optimality gap. Thus, in the regression context, $\mathbb E\|\hat\beta-\beta\|^2$ is the mean squared error (MSE) of the estimator $\hat\beta$, while in optimization this is the expected optimality gap. It is shown \citep{li1986asymptotic,li1987asymptotic} that using a regularizing penalty $\|\beta\|^2$ can improve the MSE of $\hat\beta$. However, Theorem \ref{regularization thm} shows it would not improve $\|\hat\beta-\beta\|^2$ asymptotically.

\subsection{Distance-Based DRO}\label{sec:DRO} 
Consider again that $\psi(x,P)=E_P[h(x,\xi)]$. The DRO approach entails that, under uncertainty on $P$, we obtain solution $\hat x_n^{DRO,\lambda}$ from solving
\begin{equation}
\min_{x\in\mathcal X}\max_{Q\in\mathcal U_\lambda}E_Q[h(x,\xi)]\label{DRO}
\end{equation}
where $\mathcal U_\lambda$ is a so-called uncertainty set or ambiguity set on the space of probability distributions which, at least intuitively, is believed to contain the ground-truth $P$ with high likelihood. The parameter $\lambda$ signifies the size of the set. By properly choosing $\lambda$ and accounting for the worst-case scenario in the inner maximization, \eqref{DRO} is hoped to output a higher-quality solution. DRO is gaining surging popularity in recent years. It can be viewed as a generalization of the classical deterministic RO \citep{ben2009robust, bertsimas2011theory} where the uncertain parameter in an optimization problem is now the underlying probability distribution in a stochastic problem.

Before we present our impossibility result regarding DRO over EO, let us first frame the landscape of DRO and reason the DRO types that have a legitimate possibility of beating EO statistically. In the DRO literature, the choice of $\mathcal U_\lambda$ can be categorized roughly into two groups. The first group is based on partial distributional information, such as moment and support \citep{ghaoui2003worst,delage2010distributionally,goh2010distributionally,wiesemann2014distributionally,hanasusanto2015distributionally}, shape \citep{popescu2005semidefinite,van2016generalized,li2016ambiguous,lam2017tail,chen2020discrete} and marginal distribution \citep{chen2018distributionally,doan2015robustness,dhara2021worst}. This approach has proven useful in robustifying decisions when facing limited distributional information, or when data is scarce, e.g., in the extremal region. In such cases, by using $\mathcal U_\lambda$ that captures the known partial information, the DRO guarantees a worst-case performance bound on $E_P[h(\hat x_n^{DRO,\lambda},\xi)]$ (by using the outer objective value, namely$E_P[h(\hat x_n^{DRO,\lambda},\xi)]\leq\max_{Q\in\mathcal U_\lambda}E_Q[h(\hat x_n^{DRO,\lambda},\xi)]$). Additionally, if $\mathcal U_\lambda$ is calibrated from data to be a high-confidence region in containing $P$, then such a worst-case bound holds with at least the same statistical confidence level (e.g., \citealt{delage2010distributionally}). However, from a large-sample standpoint, the set $\mathcal U_\lambda$ constructed in these approaches typically bears intrinsic looseness due to the use of only partial distributional information, and consequently the obtained solution $\hat x_n^{DRO,\lambda}$ does not converge to $x^*$ as $n$ grows (regardless of how we choose $\lambda$).

The second group of $\mathcal U_\lambda$ comprises neighborhood balls in the probability space, namely $\mathcal U_\lambda=\{Q:D(Q,P^0)\leq\lambda\}$ for some statistical distance $D(\cdot,\cdot)$ between two probability distributions, baseline distribution $P^0$, and neighborhood size $\lambda>0$. Common choices of $D$ include the $\phi$-divergence class \citep{ben2013robust,bertsimas2018robust,bayraksan2015data,jiang2016data,lam2016robust,lam2018sensitivity}, and the Wasserstein distance \citep{gao2016distributionally,chen2018robust, esfahani2018data,blanchet2019quantifying}. When the ball center $P^0$ and the ball size $\lambda$ are chosen properly in relation to the data size, it can be guaranteed that a worst-case bound holds for  $E_P[h(\hat x_n^{DRO},\xi)]$ with high confidence and, moreover, $\hat x_n^{DRO}$ converges to $x^*$. In other words, such DRO can provide statistically \emph{consistent} solutions. For this reason, in the following we will study this approach. In particular, we will present $\phi$-divergence-based DRO in detail, and then connect Wasserstein-based DRO to the result in Section \ref{sec:regularization}.



We consider $D$ represented by a $\phi$-divergence and the ball center $P^0$ taken as the empirical distribution $\hat P_n$. For any two distributions $Q$ and $Q'$ on the same domain, let $L=dQ/dQ'$ be the Radon-Nikodym derivative or the likelihood ratio between $Q$ and $Q'$. Then $D$ is defined by
\begin{equation}
D(Q,Q')=E_{Q'}[\phi(L)]\label{empirical div}
\end{equation}
where $\phi:\mathbb R_+\to\mathbb R_+$ is a function with certain properties (described momentarily). For example, when $\phi(x)=x\log x-x+1$, \eqref{empirical div} gives the Kullback-Leibler (KL) divergence, and $\phi(x)=(x-1)^2$ gives the $\chi^2$-distance. Note that such a $D$ is only well-defined when $Q$ is absolutely continuous with respect to $Q'$. Thus when we set the ball center as $\hat P_n$, the use of $D(Q,\hat P_n)$ is valid only either in the discrete-support case \citep{ben2013robust,bayraksan2015data} or functions as an ``empirical" $\phi$-divergence that confines distribution $Q$ to the support generated by data \citep{lam2017empirical,lam2019empirical,duchi2021}. In the first case, the resulting set $\mathcal U_\lambda$ can be made to contain $P$ with high confidence when $\lambda$ is chosen large enough. However, it is shown that such a calibration approach can be conservative, in the sense that a smaller $\lambda$ is sufficient for guaranteeing an upper bound for $E_P[h(\hat x_n^{DRO,\lambda},\xi)]$ (i.e., there is a looseness when translating the confidence guarantee for the uncertainty set $\mathcal U_\lambda$ in covering $P$, to the confidence level for bounding $E_P[h(\hat x_n^{DRO,\lambda},\xi)]$). The empirical $\phi$-divergence reduces this looseness by looking at the asymptotic expansion of $\max_{Q\in\mathcal U_\lambda}E_P[h(x,\xi)]$ in terms of $\lambda$ \citep{lam2016robust,gotoh2018robust,duchi2019variance} and connecting to the empirical likelihood theory \citep{wang2016likelihood,lam2017empirical,lam2019empirical,duchi2021}. In this way, even though the uncertainty set may never contain the true distribution (e.g., when the true distribution is continuous), the resulting DRO can still provide a bounding guarantee on $E_P[h(\hat x_n^{DRO,\lambda},\xi)]$ with high confidence. Finally, we note that there are other proposals in the literature that construct $P^0$ via density estimation and calibrate $\lambda$ by estimating the divergence between data and the constructed $P^0$ (e.g., \citealt{wang2009divergence,poczos2012nonparametric}). However, while being consistent, the convergence rate is slow and as a result this approach gives a very conservative $\mathcal U_\lambda$ \citep{hong2020learning} and in turn solution $\hat x_n^{DRO,\lambda}$.

To proceed, we introduce conditions on the function $\phi:\mathbb R_+\to\mathbb R_+$. 

\begin{assumption}[Smoothness of $\phi$-function]
$\phi(t)$ is convex for $t\geq0$, $\phi(1)=0$, and $\phi(t)$ is twice continuously differentiable at $t=1$ with $\phi''(1)>0$. \label{phi assumption}
\end{assumption}

Let $\phi^*(s)=\sup_{t\geq0}\{st-\phi(t)\}$ be the convex conjugate of $\phi$. Thanks to the properties of $\phi$ above, standard convex analysis \citep{rockafellar2015convex} gives the following:
\begin{lemma}[Smoothness of conjugate $\phi$-function]
Under Assumption \ref{phi assumption}, $\phi^*$ is twice continuously differentiable with $\phi^*(0)=0$, ${\phi^*}'(0)=1$ and ${\phi^*}''(0)=1/\phi''({\phi^*}'(0))=1/\phi''(1)>0$. \label{phi lemma}
\end{lemma}

Together with the regularity conditions on the optimization problem in Appendix \ref{sec:DRO proof}, we have the following result:
\begin{theorem}[Impossibility of improving EO via divergence DRO]
Consider $Z(x)=E_P[h(x,\xi)]$ with a minimizer $x^*$. Suppose we collect i.i.d. observations $\xi_1,\ldots,\xi_n$ with empirical distribution $\hat P_n$, and let $\hat x_n^{EO}$ be a minimizer of $E_{\hat P_n}[h(x,\xi)]$ and $\hat x_n^{D-DRO,\lambda}$ a minimizer of \eqref{DRO} with $\mathcal U_\lambda=\{Q:D(Q,\hat P_n)\leq\lambda\}$ and $D$ defined in \eqref{empirical div}. Under Assumption \ref{phi assumption} and Assumptions \ref{KKT}--\ref{Donsker DRO} in Appendix \ref{sec:DRO proof}, we have
\begin{eqnarray}
&&\hat x_n^{D-DRO,\lambda_n}-x^*\notag\\
&=&-(E_P[\nabla^2h(x^*,\xi)])^{-1}\langle\nabla h(x^*,\xi),\hat P_n-P\rangle-\sqrt{\lambda_n{\phi^*}''(0)}(E_P[\nabla^2h(x^*,\xi)])^{-1}\frac{Cov_P(h(x^*,\xi),\nabla h(x^*,\xi))}{\sqrt{Var_P(h(x^*,\xi))}}{}\notag\\
&&{}+o_p\left(\frac{1}{\sqrt n}+\sqrt\lambda_n\right)\label{DRO main result}
\end{eqnarray}
as $n\to\infty$ and $\lambda_n\to0$ (at any rate relative to $n$). Hence, if in addition Assumptions \ref{optimality} and \ref{nondegenerate} hold with $K=-\sqrt{{\phi^*}''(0)}(E_P[\nabla^2h(x^*,\xi)])^{-1}Cov_P(h(x^*,\xi),\nabla h(x^*,\xi))/\sqrt{Var_P(h(x^*,\xi))}$, then $n\mathcal G(\hat x_n^{EO})$ is asymptotically second-order stochastically dominated by $n\mathcal G(\hat x_n^{D-DRO,\lambda_n})$. 
\label{thm:main DRO}
\end{theorem}

In Theorem \ref{thm:main DRO}, note that $\sqrt{\lambda_n}$ plays the role of $\lambda_n$ in previous sections because of the scaling of divergence $D$. Moreover, it is known that $\max_{Q\in\mathcal U_\lambda}E_Q[h(x,\xi)]$ admits an asymptotic expansion in the form
\begin{equation}
\max_{Q\in\mathcal U_\lambda}E_Q[h(x,\xi)]\approx E_{\hat P_n}[h(x,\xi)]+\sqrt{2\lambda{\phi^*}''(0)Var_{\hat P_n}(h(x,\xi))}+\cdots\label{DRO expansion}
\end{equation}
as $\lambda\to0$ \citep{lam2016robust,dupuis2016path,gotoh2018robust,lam2018sensitivity,duchi2019variance,duchi2021}. The relation \eqref{DRO main result} can thus be obtained formally by turning the minimization of \eqref{DRO expansion} into a root-finding (or $M$-estimation; \citealt{van2000asymptotic} Chapter 5) problem via the first-order optimality condition, and then applying the delta method. Nonetheless, the precise development needs more technicality, with the detailed proof of Theorem \ref{thm:main DRO} in Appendix \ref{sec:DRO proof}. We also point out that \eqref{DRO main result} can be viewed as a generalization of \cite{duchi2019variance} that considers the special case of $\chi^2$-distance for which there is no higher-order term in \eqref{DRO expansion}, and thus our proof directly uses the Karush–Kuhn–Tucker (KKT) condition instead of starting with the Taylor expansion as in \cite{duchi2019variance}. Moreover, \eqref{DRO expansion} also relates to \cite{gotoh2021calibration} that considers the expectation and variance of the cost function of the obtained solution (see Section \ref{sec:bias-var}). 

Next we consider $D$ as the Wasserstein distance. The $p$-th order Wasserstein or optimal transport distance is defined by
\begin{equation}
D(Q,Q')=\inf_{\pi\in\Pi(Q,Q')}(E_{\pi}||\xi-\xi'\|^p)^{1/p}\label{Wasserstein def}
\end{equation}
where $\Pi(Q,Q')$ denotes the set of all distributions with marginals $Q$ and $Q'$, and $(\xi,\xi')$ is distributed according to $\pi$. The definition \eqref{Wasserstein def} has a dual representation in terms of integral probability metric \citep{sriperumbudur2012empirical}, and the norm $\|\cdot\|$ there can be replaced with more general transportation cost functions (e.g., \citealt{blanchet2019quantifying}). The rich structural properties of Wasserstein DRO has facilitated its tight connection with machine learning and statistics \citep{kuhn2019wasserstein,rahimian2019distributionally,blanchet2019profile,gao2017wasserstein,shafieezadeh2019regularization,chen2018robust}.

Here we focus on the 1-st order Wasserstein distance ($p=1$), and set the ball center $P^0$ as $\hat P_n$, so that $\mathcal U_\lambda=\{Q:D(Q,\hat P_n)\leq\lambda\}$ where $D$ is defined in \eqref{Wasserstein def}. In this case, it is known that, when $\Xi$ is in the real space and $h(x,\xi)$ is convex in $\xi$, the worst-case expectation in \eqref{DRO} possesses a Lipschitz-regularized reformulation as
\begin{equation}
\max_{Q\in\mathcal U_\lambda}E_Q[h(x,\xi)]=E_{\hat P_n}[h(x,\xi)]+\lambda\text{Lip}(h(x,\cdot))\label{regularization Wasserstein}
\end{equation}
for any $x$, where $\text{Lip}(h(x,\cdot)$ is the Lipschitz modulus of $h(x,\cdot)$ given by $\text{Lip}(h(x,\cdot)=\sup_{\xi\neq\xi'}|h(x,\xi)-h(x,\xi')|/||\xi-\xi'\|$. From \eqref{regularization Wasserstein}, we immediately obtain the following:
\begin{corollary}[Impossibility of improving EO via Wasserstein DRO]
Consider $Z(x)=E_P[h(x,\xi)]$ with a minimizer $x^*$, $h(x,\cdot)$ is convex for any $x\in\mathcal X$, and $\xi\in\Xi$ in the real space. Suppose we collect i.i.d. observations $\xi_1,\ldots,\xi_n$ with empirical distribution $\hat P_n$, and let $\hat x_n^{EO}$ be a minimizer of $E_{\hat P_n}[h(x,\xi)]$ and  $\hat x_n^{W-DRO,\lambda}$ a minimizer of \eqref{DRO} with $\mathcal U_\lambda=\{Q:D(Q,\hat P_n)\leq\lambda\}$ and $D$ defined in \eqref{Wasserstein def} with $p=1$. Under Assumptions \ref{regularity loss}--\ref{first-order assumption} in Appendix \ref{sec:proof regularization} where we set $R(x)=\text{Lip}(h(x,\cdot))$, we have
\begin{eqnarray*}
&&\hat x_n^{W-DRO,\lambda_n}-x^*\\
&=&-\langle (E_P[\nabla^2h(x^*,\xi)])^{-1}\nabla h(x^*,\xi),\hat P_n-P\rangle-\lambda_n(E_P[\nabla^2h(x^*,\xi)])^{-1}\nabla \text{Lip}(h(x^*,\cdot)){}\\
&&{}+o_p\left(\frac{1}{\sqrt n}+\lambda_n\|\nabla \text{Lip}(h(x^*,\cdot))\|\right)
\end{eqnarray*}
as $n\to\infty$ and $\lambda_n\to0$ (at any rate relative to $n$). Hence, if in addition Assumptions \ref{optimality} and \ref{nondegenerate} hold with $K=-(E_P[\nabla^2h(x^*,\xi)])^{-1}\nabla \text{Lip}(h(x^*,\cdot))$, then $n\mathcal G(\hat x_n^{EO})$ is asymptotically second-order stochastically dominated by $n\mathcal G(\hat x_n^{W-DRO,\lambda_n})$. 
\label{Wasserstein thm}
\end{corollary}

Corollary \ref{Wasserstein thm} is an immediate consequence of Theorem \ref{regularization thm} by observing the Lipschitz-regularized reformulation of Wasserstein DRO. Though we do not pursue further in obtaining sufficient conditions for the required regularity properties of $\text{Lip}(h(x,\cdot))$, we point out that this can  be written explicitly in terms of the dual norm of $\|\cdot\|$ in \eqref{Wasserstein def} for interesting examples such as regression and classification problems \citep{esfahani2018data,blanchet2019profile,kuhn2019wasserstein}. Moreover, more generally (i.e., $p\neq1$ and $h(x,\cdot)$ not necessarily convex), the worst-case objective $\max_{Q\in\mathcal U_\lambda}E_Q[h(x,\xi)]$ admits an expansion similar to \eqref{DRO expansion} where the first-order term is $\lambda V(h(x,\cdot))$, with $V(h(x,\cdot))$ being some variability measure of $h$ in which $\text{Lip}(h(x,\cdot))$ is a special case \citep{gao2017wasserstein}.

\subsection{Parametric Optimization}\label{sec:parametric}
We have focused on nonparametric problems thus far in the previous subsections. Our framework applies similarly if we confine the unknown distribution $P$ to a parametric model, and we will present an impossibility result in improving parametric optimization using regularization. 

In this parametric case, we write $Z(x)=\psi(x,\theta^*)$ where $\theta^*$ is an unknown parameter in the model, e.g., $\psi(x,\theta)=E_\theta[h(x,\xi)]$ where $E_\theta[\cdot]$ denotes the expectation under a parametric distribution, say $P_\theta$, that generates $\xi$ and $\theta^*$ is the true value of $\theta$. Evidently, we can write $\theta^*=\theta(P)$ where $\theta(\cdot)$ is viewed as a function on the true distribution $P$ (e.g., $\theta(\cdot)$ can be the solution to the score function equation), so we can define $\tilde\psi(x,Q)=\psi(x,\theta(Q))$ and thus $Z(x)=\tilde\psi(x,P)$. Suppose, given data, we use a consistent estimator $\hat\theta_n$ to plug into $\psi(x,\hat\theta_n)$, then we can write $\hat\theta_n=\theta(\hat P_n)$. Thus we have $\min_{x\in\mathcal X}\tilde\psi(x,\hat P_n)$ as our EO procedure, which reduces to the same setting as the previous subsections. Note that the above discussion implicitly assumes the model is correctly specified, so that the estimator $\hat\theta_n$ converges to the true parameter value $\theta^*$. Nonetheless, the framework still works even if the model is misspecified, in which case $\theta^*$ denotes the minimizer of the Kullback-Leibler divergence between the true model and the considered model class, and $\hat\theta_n$ is then consistent in estimating $\theta^*$ under mild conditions.

To make our discussion precise, let $x^*$ be an optimal solution to $\min_{x\in\mathcal X}\{Z(x)=\psi(x,\theta^*)\}$ in this parametric setting, where $\theta^*\in\mathbb R^m$ is the true parameter (either when the model is correctly specified, or interpreted as the minimizer of the Kullback-Leibler divergence when it is incorrectly specified). To avoid ambiguity, we write $\nabla_x\in\mathbb R^d$ and $\nabla_x^2\in\mathbb R^d$ as the gradient and Hessian taken with respect to $x$, and $\nabla_{x\theta}\in\mathbb R^{d\times m}$ as the cross-gradients with respect to $\theta$ and $x$. We write $\hat x_n^{P-EO}\in\text{argmin}_{x\in\mathcal X}\psi(x,\hat\theta_n)$ where $\hat\theta_n$ is obtained from an estimation procedure on the probability distribution model that satisfies the asymptotic behavior in \eqref{parametric CLT} depicted momentarily. We also consider a penalty function $R:\mathbb R^d\to\mathbb R$, and the regularized formulation
$$\hat x_n^{P-Reg,\lambda}\in\text{argmin}_{x\in\mathcal X}\{\psi(x,\hat\theta_n)+\lambda R(x)\}$$
We have the following result:
\begin{theorem}[Impossibility of improvement over EO in parametric optimization]
Suppose Assumptions \ref{assumption:parametric objective} and \ref{assumption:penalty} in Appendix \ref{sec:parametric proof} hold. Also suppose $\hat\theta_n$ is obtained from data such that
\begin{equation}\hat\theta_n-\theta^*=\langle IF_\theta,\hat P_n-P\rangle+o_p\left(\frac{1}{\sqrt n}\right)\label{parametric CLT}
\end{equation}
where $IF_\theta(\xi):\Xi\to\mathbb R^d$ has $Cov_P(IF_\theta(\xi))$ that is finite. Then
$$\hat x_n^{P-Reg,\lambda}-x^*=-\langle(\nabla_x^2\psi(x^*,\theta^*))^{-1}\nabla_{x\theta}\psi(x^*,\theta^*)IF_\theta,\hat P_n-P\rangle-\lambda(\nabla_x^2\psi(x^*,\theta^*))^{-1}\nabla_xR(x)+o_p\left(\frac{1}{\sqrt n}+\lambda\right)$$
as $n\to\infty$ and $\lambda\to0$. Hence, if in addition Assumptions \ref{optimality} and \ref{nondegenerate} hold with $K=-(\nabla_x^2\psi(x^*,\theta^*))^{-1}\nabla_xR(x)$, then $n\mathcal G(\hat x_n^{P-EO})$ is asymptotically second-order stochastically dominated by $n\mathcal G(\hat x_n^{P-Reg,\lambda_n})$ for any $\lambda_n\to0$. \label{thm:parametric}
\end{theorem}

The asymptotic \eqref{parametric CLT} can be ensured by standard conditions. For instance, in the case where there is no model misspecification and we use maximum likelihood estimator for $\hat\theta_n$, \eqref{parametric CLT} holds under the Lipschitzness of the log-likelihood function, with $IF_\theta(\cdot)=I_{\theta^*}^{-1}s_{\theta^*}(\cdot)$ where $s_{\theta^*}(\cdot)$ is the score function and $I_{\theta^*}$ is the Fisher information matrix (see, e.g., \citealt{van2000asymptotic} Theorem 5.39).

\subsection{Bayesian Optimization}\label{sec:Bayesian}
We generalize the discussion in Section \ref{sec:parametric} to Bayesian optimization, in the following sense. We impose a prior distribution on the unknown parameter $\theta$, and obtain $\hat x_n^{P-Bay,\lambda}$ from
\begin{equation}
\hat x_n^{P-Bay,\lambda}\in\text{argmin}_{x\in\mathcal X}\{E_{\Theta|D_n}[\psi(x,\Theta)]+\lambda R(x)\}\label{Bayesian solution}
\end{equation}
where $E_{\Theta|D_n}[\cdot]$ is the posterior distribution of $\theta$ given the collection of data $D_n=\{\xi_1,\ldots,\xi_n\}$, and we denote $\Theta$ as the random variable distributed under this posterior distribution. In other words, $\hat x_n^{P-Bay,0}$ optimizes the posterior expectation of the original objective function $\psi(x,\theta)$, while $\hat x_n^{P-Bay,\lambda}$ imposes additionally a regularizing penalty $R$ with the regularization parameter $\lambda$.
We have the following result:
\begin{theorem}[Impossibility of improvement over EO via Bayesian optimization]
Suppose all assumptions in Theorem \ref{thm:parametric} and Assumptions \ref{assumption:Bayesian}, \ref{assumption:Bayesian penalty} and \ref{assumption:Bayesian solution} in Appendix \ref{sec:Bayesian proof} hold. In addition, assume that
\begin{equation}
\left\|P_{\Theta|D_n}-N\left(\hat\theta_n,\frac{J}{n}\right)\right\|_{TV}\stackrel{p}{\to}0\label{posterior convergence}
\end{equation}
for some $\hat\theta_n$ where $\|\cdot\|_{TV}$ denotes the total variation distance, $P_{\Theta|D_n}$ is the posterior distribution of $\theta$, and $J$ is some covariance matrix. 
Moreover, assume that $\sqrt n(\Theta-\hat\theta_n)|D_n$ is uniformly integrable a.s. as $n\to\infty$. 
Then
$$\hat x_n^{P-Bay,\lambda}-x^*=-\langle(\nabla_x^2\psi(x^*,\theta^*))^{-1}\nabla_{x\theta}\psi(x^*,\theta^*)IF_\theta,\hat P_n-P\rangle-\lambda(\nabla_x^2\psi(x^*,\theta^*))^{-1}\nabla_xR(x)+o_p\left(\frac{1}{\sqrt n}+\lambda\right)$$
as $n\to\infty$ and $\lambda\to0$. Hence, if in addition Assumptions \ref{optimality} and \ref{nondegenerate} hold with $K=-(\nabla_x^2\psi(x^*,\theta^*))^{-1}\nabla_xR(x)$, then $n\mathcal G(\hat x_n^{P-EO})$ is asymptotically second-order stochastically dominated by $n\mathcal G(\hat x_n^{P-Bay,\lambda_n})$ for any $\lambda_n\to0$. \label{thm:Bayesian}
\end{theorem}

Note that the asymptotic expressions of $\hat x_n^{P-Bay,\lambda}-x^*$ in Theorem \ref{thm:Bayesian} and $\hat x_n^{P-Reg,\lambda}-x^*$ in Theorem \ref{thm:parametric} are the same. This in particular implies that the simple Bayesian solution, $\hat x_n^{P-Bay,0}$, performs asymptotically equivalently to the EO solution and, moreover, $\hat x_n^{P-Bay,\lambda}$ and $\hat x_n^{P-Reg,\lambda}$ perform asymptotically equivalently.

Finally, we mention that the convergence \eqref{posterior convergence} is standard, as guaranteed by the Berstein-von Mises Theorem (e.g., Theorem 10.1 and the discussion on P.144 in \citealt{van2000asymptotic}).

\section{Comparisons with Existing Literature}\label{sec:previous}
As we have mentioned in the introduction, in the literature data-driven formulations other than EO have been proposed for good reasons. A focal point in recent years has been on DRO, which have been studied by a range of interesting works. Below we compare our framework with these DRO works in Sections \ref{sec:abs} and \ref{sec:bias-var}, but noting that our framework in Section \ref{sec:main} applies more generally to all ``expanded" procedures in addition to DRO. In Sections \ref{sec:cr} and \ref{sec:op} we connect our result to semiparametric inference in statistics, and discuss alternate settings in which other approaches than EO could offer advantages. 

\subsection{Absolute Performance Bounds and Out-of-Sample Disappointment}\label{sec:abs}
A common guarantee studied in the DRO literature is that the worst-case or robust objective value provides an upper bound on the true objective value evaluated at the obtained solution. Consider the expected value minimization problem $\min_{x\in\mathcal X}\{Z(x)=E_P[h(x,\xi)]\}$. Suppose $\mathcal U_\lambda$ is an uncertainty set that contains the true distribution $P$. Then the DRO solution $\hat x^{DRO}$ obtained from solving $\min_{x\in\mathcal X}\max_{Q\in\mathcal U_\lambda}E_Q[h(x,\xi)]$ must satisfy
$$Z(\hat x^{DRO})\leq\max_{Q\in\mathcal U_\lambda}E_Q[h(\hat x^{DRO},\xi)]$$
Moreover, when data are observed, suppose $\mathcal U_\lambda$ is constructed as a high-confidence region for $P$, i.e., $\mathbb P(P\in\mathcal U_\lambda)\geq1-\alpha$ for some confidence level $1-\alpha$. Then this confidence level can be translated into at least the same confidence on a bound of the performance of $\hat x^{DRO}$ via the worst-case objective value, given by
\begin{equation}
\mathbb P\left(Z(\hat x^{DRO})\leq\max_{Q\in\mathcal U_\lambda}E_Q[h(\hat x^{DRO},\xi)]\right)\geq1-\alpha\label{absolute bound}
\end{equation}
where $\mathbb P$ denotes the probability with respect to the data. The above argument can be finite-sample or asymptotically argued (in the latter case, the asymptotic confidence guarantee of the uncertainty set would translate to the asymptotic confidence guarantee of the performance bound). Such type of guarantees has been studied in, e.g., \cite{delage2010distributionally,goh2010distributionally,hanasusanto2015distributionally,lam2017tail,ben2013robust,bertsimas2018robust,jiang2016data,esfahani2018data}. \cite{van2020data} and \cite{sutter2020general} in particular call the complement of the probability in the left hand side of \eqref{absolute bound} the out-of-sample disappointment.

Our results in Sections \ref{sec:main} and \ref{sec:DRO} differ from the bound \eqref{absolute bound} in two important aspects. First is that we are measuring the quality of solution by a ranking of the true objective value. That is, an obtained solution $\hat x$ that has a smaller value of $Z(\hat x)$ is regarded as more desirable. This is different from \eqref{absolute bound} that guarantees the validity of the \emph{estimated} objective value in bounding the true value of an obtained solution. Note that the latter validity does not necessarily imply the obtained solution performs better in the \emph{true} objective value. In fact, we show that DRO cannot be superior to EO in the latter aspect, at least in the large-sample regime that we consider.

Regarding which criterion, \eqref{absolute bound} or ours, should an optimizer use, it may depend on the particular situation of interest. In terms of comparing solution quality, we believe there should be little argument against our criterion of ranking the attained true objective value, as this appears the most direct measurement of solution performance. On the other hand, in some high-stake situations, an optimizer may want to obtain a reliable upper estimate, or to ensure a low enough upper bound, of the attained objective value, in which case the conventional DRO guarantee \eqref{absolute bound} would be useful. 

Our second distinction from bound \eqref{absolute bound} is that we study the error \emph{relative} to the oracle best solution, i.e., our approximation is on the optimality gap $Z(\hat x)-Z(x^*)$ for an obtained solution $\hat x$. A claimed drawback of using bound \eqref{absolute bound} is that it could be loose, thus unable to detect the over-conservativeness of DRO (e.g., one can simply take $\mathcal U_\lambda$ to be extremely large, so that \eqref{absolute bound} trivially holds). This latter criticism is resolved to an extent by a series of work that shows that, by choosing $\mathcal U_\lambda$ properly, the worst-case objective value $\max_{Q\in\mathcal U_\lambda}E[h(x,\xi)]$ differs from the true value $E_P[h(x,\xi)]$ only by a small amount. This includes the empirical likelihood \citep{lam2017empirical,lam2019empirical,duchi2019variance,duchi2021}, Bayesian \citep{gupta2015}, and large deviations perspective \citep{van2020data} in divergence DRO, and the profile likelihood and variability regularization for Wasserstein DRO \citep{blanchet2019profile,gao2017wasserstein1}. Moreover, it can be shown that certain divergence DRO gives rise to the best possible bound in the form of \eqref{absolute bound}, which is argued via a match of the statistical performance with the CLT \citep{lam2019empirical,gupta2015,duchi2021} or the Bernstein bound \citep{duchi2019variance}, or a ``meta-optimization" that gives the tightest such bound subject to an allowed large deviations rate \citep{van2020data}. Nonetheless, all these works mainly focus on an upper bound on $Z(\hat x)$ instead of the optimality gap $Z(\hat x)-Z(x^*)$, the latter being more challenging due to the unknown oracle true optimal value $Z(x^*)$. 

\subsection{Bias-Variance Tradeoff}\label{sec:bias-var}
Consider the expected value minimization problem $\min_{x \in\mathcal X}\{Z(x)=E_P[h(x,\xi)]\}$ where $h$ is some loss function. \cite{gotoh2018robust,gotoh2021calibration} studies a regularized formulation using $\phi$-divergence as the penalty, defined as (in the notation of the current paper):
\begin{equation}
\min_{x\in\mathcal X}\max_Q\left\{E_Q[h(x,\xi)]+\frac{1}{\lambda}D(Q,\hat P_n)\right\}\label{gotoh}
\end{equation}
where $Q$ lies in the relevant space of probability distributions. This formulation is akin to the divergence DRO we discussed in Section \ref{sec:DRO}, but using the Lagrangian formulation directly instead of starting with the notion of uncertainty set. The parameter $\lambda$ in \eqref{gotoh} plays a similar role as the ball size of our uncertainty set.

The main insight from \cite{gotoh2018robust,gotoh2021calibration} is a desirable improvement on the bias-variance tradeoff using the solution obtained in \eqref{gotoh}, which we call a ``Lagrangian (L)-DRO" solution $\hat x_n^{L-DRO,\lambda}$. More precisely, \cite{gotoh2018robust} concludes (Theorem 5.1 therein) that
\begin{eqnarray}
&&E_{\hat P_n}[h(\hat x_n^{L-DRO,\lambda},\xi)]\notag\\
&=&E_{\hat P_n}[h(\hat x_n^{EO},\xi)]+\frac{\lambda^2}{2{\phi^*}''(0)^2}Cov_{\hat P_n}(h(\hat x_n^{EO},\xi),\nabla h(\hat x_n^{EO},\xi))^\top(E_{\hat P_n}[\nabla^2h(\hat x_n^{EO},\xi)])^{-1}Cov_{\hat P_n}(h(\hat x_n^{EO},\xi),\nabla h(\hat x_n^{EO},\xi)){}\notag\\
&&{}+o(\lambda^2)\label{bias L-DRO}
\end{eqnarray}
and
\begin{eqnarray}
&&Var_{\hat P_n}[h(\hat x_n^{L-DRO,\lambda},\xi)]\notag\\
&=&Var_{\hat P_n}[h(\hat x_n^{EO},\xi)]-\frac{2\lambda}{{\phi^*}''(0)}Cov_{\hat P_n}(h(\hat x_n^{EO},\xi),\nabla h(\hat x_n^{EO},\xi))^\top(E_{\hat P_n}[\nabla^2h(\hat x_n^{EO},\xi)])^{-1}Cov_{\hat P_n}(h(\hat x_n^{EO},\xi),\nabla h(\hat x_n^{EO},\xi)){}\notag\\
&&{}+o(\lambda)\label{var L-DRO}
\end{eqnarray}
as $\lambda\to0$, and \cite{gotoh2021calibration} develops these expansions further to center at $E_P[h(\hat x_n^{EO},\xi)]$ and $Var_P(h(\hat x_n^{EO},\xi))$. From \eqref{bias L-DRO} and \eqref{var L-DRO}, we see that while $L-DRO$ deteriorates the expected loss, it reduces the variance of the loss by a large magnitude ($\lambda$ versus $\lambda^2$), giving an overall improvement on the mean squared error. In this sense, injecting a small $\lambda>0$ in $L-DRO$ is desirable compared to EO with $\lambda=0$.

Putting aside the technical differences in using L-DRO versus the common DRO in \eqref{DRO} (the latter requires an extra layer of analysis in the Lagrangian reformulation), we point out two conceptual distinctions between \cite{gotoh2018robust,gotoh2021calibration} and our results in Section \ref{sec:main} and \ref{sec:DRO}. First is the criterion in measuring the quality of an obtained solution $\hat x$, in particular the role of the variance of the loss function $h(x,\xi)$, $Var_P(h(x,\xi))$. Our criterion is in terms of the achieved optimality gap, or equivalently $Z(\hat x)$, where $Z(\cdot)$ is the \emph{true} objective function of the original optimization. In particular, any risk-aware consideration should be already incorporated into the construction of the loss $h$. When an obtained solution $\hat x$ is used in many future test cases, an estimate of $Z(\hat x)$, using $n_{test}$ test data points, has a variance given by $Var_P(h(\hat x,\xi))/n_{test}$ (instead of $Var_P(h(\hat x,\xi))$), and thus the variance of $h$ plays a relatively negligible role. This is different from \cite{gotoh2018robust,gotoh2021calibration} who takes an alternate view that puts more weight on the variability of the loss function.

Our second main distinction from \cite{gotoh2018robust,gotoh2021calibration} is our consideration of the attained true objective value or generalization performance $Z(\hat x)$ as a random variable, and we study its risk profile by assessing its \emph{entire distribution} that exhibits the second-order stochastic dominance relation put forth in Section \ref{sec:main}. Here, the randomness of $Z(\hat x)$ comes from the statistical noise from data used in constructing the obtained solution $\hat x$. In contrast, \cite{gotoh2018robust} studies the mean of the attained objective value (and variance in the sense described above). In this latter setting, as we have seen in \eqref{expected gap}, any expanded procedure with parameter $\lambda$ satisfying Assumption \ref{solution expansion} (not only divergence DRO) would lead to a deterioration of order $\lambda^2$ from EO. This does show an inferiority of the expanded procedure, but it does not give a complete picture as the attained objective value can be better in other distributional aspects. Our Theorem \ref{main thm} stipulates that even considering the entire distribution, an expanded procedure like DRO still cannot outperform EO.

\subsection{Connections to Local Asymptotic Minimax Theorem and Cramer-Rao Bound}\label{sec:cr}
Our main result connects to the classical asymptotic minimax theorem  in semiparametric inference and the well-known Cramer-Rao bound. Consider the optimal solution $x^*$ as an estimation target, and we construct an estimator $\hat x_n$ from data. It is known from the local asymptotic minimax theorem (\citealt{kosorok2007introduction} Theorem 18.4) that, for any so-called regular estimator $\hat x_n$ (\citealt{kosorok2007introduction} Section 18.1), $\sqrt n(\hat x_n-x^*)$ asymptotically dominates $N(0,\Sigma)$ where $\Sigma=E_P[\nabla^2h(x^*,\xi)])^{-1}Cov_P(\nabla h(x^*,\xi))E_P[\nabla^2h(x^*,\xi)])^{-1}$. Note that the latter variable is precisely the limit of $\sqrt n(\hat x_n^{EO}-x^*)$ ($\Sigma/n$ is the the variance of the first term in \eqref{regularization expression}). This asymptotic dominance is in the sense that $\mathbb E[\ell(\sqrt n(\hat x_n-x^*))]\geq\mathbb E[\ell(N(0,\Sigma))]$ under a proper asymptotic limit, for any so-called subconvex function $\ell$, namely $\ell$ such that the sublevel sets $\{y:\ell(y)\leq c\}$ for any $c\in\mathbb R$ are convex, symmetric about the origin, and closed. 

On the other hand, our second-order stochastic dominance in Theorem \ref{main thm} shows that $\mathbb E[f(n\mathcal G(\hat x_n))]$ for any considered data-driven solution $\hat x_n$ is at least $\mathbb E[f(n\mathcal G(\hat x_n^{EO})])$, for any non-decreasing convex function $f$. Noting from \eqref{new opt gap CR} that $\mathcal G(\hat x_n)\approx\frac{1}{2}(\hat x_n-x^*)^\top\nabla^2Z(x^*)(\hat x_n-x^*)$, our result implies, roughly speaking, that $\mathbb E[f(g(\sqrt n(\hat x_n-x^*))]$ is at least $\mathbb E[f(g(\sqrt n(\hat x_n^{EO}-x^*))]$ asymptotically, where $g(y)=\|(\nabla^2Z(x^*))^{1/2}y\|^2$ which is convex and thus $f\circ g$ is subconvex by the non-decreasing convex property of $f$ and the form of $g$. Our result thus, at least intuitively, coincides with the implication of the minimax theorem. However, the minimax theorem considers a perturbation of the true parameter value within a shrinking neighborhood at a specific rate  (which leads to the notion of regular estimators and also defines the plausible parameter values to be considered in the minimax regime), and it claims the dominance relation using Anderson's lemma, a result in convex geometry (\citealt{van2000asymptotic} Theorem 8.5). On the other hand, our Theorem \ref{main thm} is obtained by perturbing the data-driven procedure parametrized by $\lambda$, and analyzes the optimality gap or objective value via second-order stochastic dominance. Our main insight to claim the superiority of EO, which uses the mean-preserving spread in risk-based ranking, appears novelly beyond the semiparametric literature.







Moreover, it is perhaps revealing to cast our result in the context of the conventional Cramer-Rao bound. The latter states that, under suitable smoothness conditions, maximum likelihood is the best among all possible estimators in estimating (functions of) unknown model parameters in terms of asymptotic variance. Suppose we use a particular objective function, namely the expected log-likelihood, in our framework. The Cramer-Rao bound would intuitively conclude the expected optimality gap is best attained via EO, and thus our result can be viewed as a local generalization in that it applies to general objective functions and concludes the superiority of the entire distribution of the optimality gap when using EO.

To explain in detail, first note that \eqref{prelim expression} stipulates that the EO solution, $\hat x_n^{EO}$, satisfies
$$Z(\hat x_n^{EO})-Z(x^*)\approx\frac{1}{2}(\hat x_n^{EO}-x^*)^\top\nabla^2Z(x^*)(\hat x_n^{EO}-x^*)$$
Now consider, under suitable conditions like the ones for Section \ref{sec:examples},
\begin{align}
E[(\hat x_n^{EO}-x^*)^\top\nabla^2Z(x^*)(\hat x_n^{EO}-x^*)]&=E[\langle IF,\hat P_n-P\rangle^\top\nabla^2Z(x^*)\langle IF,\hat P_n-P\rangle]\notag\\
&\approx E\left[\frac{1}{n}\sum_{i=1}^n\nabla h(x^*,\xi_i)^\top\nabla^2Z(x^*)^{-1}\nabla^2Z(x^*)\nabla^2Z(x^*)^{-1}\frac{1}{n}\sum_{i=1}^n\nabla h(x^*,\xi_i)\right]\notag\\
&=\frac{tr(\nabla^2Z(x^*)^{-1}Cov(\nabla h(x^*,\xi)))}{n}\label{CR interim}
\end{align}
where in the approximate equality we apply Theorem \ref{regularization thm} with $\lambda=0$ to deduce $\hat x_n^{EO}-x^*\approx\nabla^2Z(x^*)^{-1}\frac{1}{n}\sum_{i=1}^n\nabla h(x^*,\xi_i)$ and noting that $\nabla^2Z(x^*)=E_P[\nabla^2h(x^*,\xi)]$, and in the last equality we use the first-order optimality condition that $E_P[\nabla h(x^*,\xi)]=0$.

Now, consider the special case where $Z(x)=E_P[h(x,\xi)]$ and  $h(x,\xi)=-\log f_x(\xi)$ for a parametric class of density $\{f_x(\cdot)\}_x$. Suppose $f_{x^*}(\xi)$ is the true density of the random variable $\xi$, whose distribution is $P$. We therefore have $Z(x)=E_P[-\log f_x(\xi)]$ and $x^*$, which is the true value of the parameter in the density function of $\xi$, satisfies $x^*=\text{argmin}_xZ(x)$ as it minimizes the Kullback-Leibler divergence from the true distribution among all densities parametrized as $f_x(\cdot)$.

Note that in this special case,
\begin{align*}
tr(\nabla^2Z(x^*)^{-1}Cov_P(\nabla h(x^*,\xi)))&=tr(\nabla^2Z(x^*)^{-1}E_P[\nabla h(x^*,\xi))^2]\\
&=tr(\nabla^2Z(x^*)^{-1}E_P[\nabla\log f_{x^*}(\xi))^2]\\
&=tr(\nabla^2Z(x^*)^{-1}(-E_P[\nabla^2\log f_{x^*}(\xi)]))\\
&=tr(\nabla^2Z(x^*)^{-1}\nabla^2Z(x^*))\\
&=d
\end{align*}
where we use the first-order optimality condition $E_P[\nabla h(x^*,\xi)]=0$ in the first equality and the equivalent expression of the Fisher information matrix in the third equality. This shows that
\eqref{CR interim} becomes $d/n$ in this special case.

Now consider $\hat x_n^{EO}$ as an estimator of $x^*$, or that $(\nabla^2Z(x^*))^{1/2}\hat x_n^{EO}$ as an estimator of $(\nabla^2Z(x^*))^{1/2}x^*$. The Cramer-Rao bound states that, for any nearly unbiased estimator $(\nabla^2Z(x^*))^{1/2}\hat x_n$, we have
$$tr(Cov_P((\nabla^2Z(x^*))^{1/2}\hat x_n))\geq tr(\nabla\psi(x^*)'(nI)^{-1}\nabla\psi(x^*))$$
where $\psi(x)=E_x[(\nabla^2Z(x^*))^{1/2}\hat x_n]=(\nabla^2Z(x^*))^{1/2}E_x\hat x_n$, with $E_x[\cdot]$ being the expectation taken under $f_x(\cdot)$, and $I$ is the Fisher information matrix given by $I=-E_P[\nabla^2\log f_{x^*}(\xi)]=\nabla^2Z(x^*)$. Now, consider any estimator $\hat x_n$ that satisfies $E_x[\hat x_n]= x+O(1/n)$ (such as all the data-driven solutions we discussed in Section \ref{sec:examples} using small enough $\lambda_n$, i.e., $\lambda_n=O(1/n)$). We have $\nabla\psi(x^*)\approx(\nabla^2Z(x^*))^{1/2}$ and 
$$E\|(\nabla^2Z(x^*))^{1/2}(\hat x_n-x^*)\|^2\approx tr(Cov_P((\nabla^2Z(x^*))^{1/2}\hat x_n))\geq {tr((\nabla^2Z(x^*))^{1/2}}^\top(nI)^{-1}(\nabla^2Z(x^*))^{1/2})=\frac{d}{n}$$
Therefore, in the case where $\xi\sim f_{x^*}(\cdot)$ and $h(x,\xi)=-\log f_x(\xi)$, the EO solution $\hat x_n^{EO}$ is approximately optimal in the sense that
the optimality gap $Z(x)-Z(x^*)$, which is approximately
$\frac{1}{2}(x-x^*)^\top\nabla^2Z(x^*)(x-x^*)$ for any estimator $x$ close to $x^*$ in expectation, has expectation $E[Z(x)-Z(x^*)]$ minimized at $\hat x_n^{EO}$ up to a negligible error. Our main result is a generalization of the above. It suggests that $\hat x_n^{EO}$ not only minimizes $Z(x)-Z(x^*)$ in expectation, but also minimizes $E[g(Z(x)-Z(x^*))]$ for any convex non-decreasing loss function $g(\cdot)$, or equivalently that $Z(\hat x_n)-Z(x^*)$  for any data-driven solution $\hat x_n$ second-order stochastically dominates $Z(\hat x_n^{EO})-Z(x^*)$. 


\subsection{Beyond Empirical Optimization}\label{sec:op}
Finally, while this paper advocates the superiority of EO in standard large-sample settings, we point out, as also mentioned in the introduction, that there are many other settings where alternate data-driven strategies than EO offer benefits. We have already mentioned some of them in Sections \ref{sec:abs} and \ref{sec:bias-var} when comparing with some existing DRO studies. We close this paper by discussing some other prominent examples and the types of procedures that are demonstrably powerful:

\emph{Finite-sample performance: }We have focused and claimed EO is best in the large-sample regime, but for small sample the situation could be different. While it may be difficult to argue there are generally superior procedures than EO (or vice versa) in finite sample, there exists documented situations where some approaches are consistently better than EO. For example, the so-called operational statistic \citep{liyanage2005practical,chu2008solving} strengthens data-driven solutions to perform better than simple choices such as EO, uniformly across all possible parameter values (in a parametric problem) and sample size. It does so by expanding the class of data-driven solutions in which a better solution is systematically searched via a second-stage optimization or Bayesian analysis. In certain structural problems, for instance those arising in inventory control, this approach provably and empirically performs better than EO.

\emph{High dimension: }We have fixed the problem dimension throughout this paper. It is well-known, however, that certain regularization schemes, e.g., $L_1$-penalty as in LASSO \citep{friedman2001elements}, are suited for high-dimensional problems that exhibit sparse structures. Through the equivalence with regularization, this also implies that certain types of Wasserstein DRO enjoy similar statistical benefits \citep{blanchet2019profile,shafieezadeh2019regularization}. 

\emph{Distortion of loss function properties: }For problems that do not satisfy our imposed smoothness conditions, the decay rate of $\mathcal G(\hat x_n^{EO})$ could be $1/\sqrt n$ instead of $1/n$. In this situation, $\chi^2$-divergence DRO, which exhibits variance regularization, can be used to boost the rate to $1/n$ in certain examples \citep{duchi2019variance}. Moreover, it is also shown that DRO using a distance induced from the reproducing kernel Hilbert space \citep{staib2019distributionally}, which relates to the so-called kernel mean matching \citep{gretton2009covariate}, could exhibit finite-sample optimality gap bounds that are free of any complexity of the loss function class, which is drastically distinct from EO's behavior \citep{zenglam21}.

\emph{Data pooling and contextual optimization: }When there are simultaneously many stochastic optimization problems to solve, it is shown that introducing a shrinkage onto EO can enhance the tradeoff between the optimality gap and instability \citep{gupta2021small,gupta2021data}. This shrinkage is in a similar spirit as the so-called James-Stein estimator \citep{cox1979theoretical} in classical statistics, in which individual estimators, or the solutions of the individual optimization problems, are adjusted by ``weighting" with a pooled estimator. Relatedly, when the considered optimization problem involves parameters or outcomes that depend on covariates, using no such covariate information, or using the naive predict-then-optimize approach, can be improved by integrating the prediction and EO steps that leads to a better ultimate objective value performance (e.g., \citealt{ban2019big,elmachtoub2021smart}).

\ACKNOWLEDGMENT{I gratefully acknowledge support from the National Science Foundation under grants CAREER CMMI-1834710 and IIS-1849280.}






\bibliographystyle{informs2014}
\bibliography{bibliography}

\ECSwitch
\ECHead{Proofs of Statements}
\section{Review of Useful Results}

The following is a classical theorem on $M$-estimation in the presence of extra parameters in the stochastic equation.
\begin{theorem}[a.k.a. Theorem 5.31 in \citealt{van2000asymptotic}]Consider, for all $x\in\mathcal X\subset\mathbb R^d$ and $\lambda\in\mathbb R^m$ such that $\|\lambda-\lambda_0\|\leq\delta$ for a given point $\lambda_0$ for some $\delta>0$ , $\varphi_{x,\lambda}(\cdot)$ is a function mapping $\Xi$ to $\mathbb R^d$. Also consider (possibly random) sequences $x_n\in\mathcal X$ and $\lambda_n$ for $n=1,2,\ldots$. Let $P$ be a distribution generating $\xi\in\Xi$ and $E_P[\cdot]$ denotes its expectation. Let $\hat P_n$ denotes the empirical distribution of i.i.d. observations $\xi_1,\ldots,\xi_n$ and $E_{\hat P_n}[\cdot]$ denotes its expectation. We assume the following conditions:
\begin{enumerate}
\item $\{\varphi_{x,\lambda}(\cdot):\|x-x^*\|\leq\delta,\|\lambda-\lambda_0\|\leq\delta\}$ is Donsker. 
\item $E_P\|\varphi_{x,\lambda}(\xi)-\varphi_{x^*,\lambda_0}(\xi)\|^2\to0$ as $(x,\lambda)\to(x^*,\lambda_0)$.
\item $E_P[\varphi_{x^*,\lambda_0}(\xi)]=0$
\item $E_P[\varphi_{x,\lambda}(\xi)]$ is differentiable at $x=x^*$, uniformly in $\lambda$ within a small neighborhood of $\lambda_0$ with non-singular derivative matrix $V_{x^*,\lambda}=\nabla E_P[\varphi_{x^*,\lambda}(\xi)]$ such that $V_{x^*,\lambda}\to V_{x^*,\lambda_0}$.
\item $\sqrt nE_{\hat P_n}[\varphi_{x_n,\lambda_n}(\xi)]=o_p(1)$
\item $(x_n,\lambda_n)\stackrel{p}{\to}(x^*,\lambda_0)$ as $n\to\infty$
\end{enumerate}
Then
$$x_n-x^*=-V_{x^*,\lambda_0}^{-1}E_P[\varphi_{x^*,\lambda_n}(\xi)]-V_{x^*,\lambda_0}^{-1}E_{\hat P_n}[\varphi_{x^*,\lambda_0}(\xi)]+o_p\left(\frac{1}{\sqrt n}+\|E_P[\varphi_{x^*,\lambda_n}(\xi)]\|\right)$$
\label{nuisance}
\end{theorem}

We also need the following classical consistency theorem in $M$-estimation.
\begin{theorem}[a.k.a. Theorem 5.9 in \citealt{van2000asymptotic}]Let $\Psi_n(\cdot)$ be a random vector-valued functions and $\Psi(\cdot)$ be a deterministic vector-valued function of $x\in\mathcal X$ such that for every $\epsilon>0$,
$$\sup_{x\in\mathcal X}\|\Psi_n(x)-\Psi(x)\|\stackrel{p}{\to}0$$
$$\inf_{x\in\mathcal X:\|x-x^*\|\geq\epsilon}\|\Psi(x)\|>0=\|\Psi(x^*)\|$$
Then any sequence $x_n$ such that $\Psi_n(x_n)=o_p(1)$ converges in probability to $x^*$.\label{consistency}
\end{theorem}

\section{Proofs for Section \ref{sec:main}}\label{sec:proof main}
\proof{Proof of Lemma \ref{explicit lemma}.}
Define
\begin{equation}
r(x)=\frac{Z(x)-Z(x^*)-\frac{1}{2}(x-x^*)^\top\nabla^2Z(x^*)(x-x^*)}{\|x-x^*\|^2}\label{interim1}
\end{equation}
for $x\neq x^*$ and $r(x^*)=0$. By Assumption \ref{optimality}, we know that $r(x)$ is continuous at $x^*$. From \eqref{interim1}, we can write
\begin{equation}
Z(\hat x_n^\lambda)-Z(x^*)=r(\hat x_n^\lambda)\|\hat x_n^\lambda-x^*\|^2+\frac{1}{2}(\hat x_n^\lambda-x^*)^\top\nabla^2Z(x^*)(\hat x_n^\lambda-x^*)\label{interim2}
\end{equation}
From Assumption \ref{solution expansion}, we have \eqref{asymptotic expansion} and \eqref{asymptotic expansion original} with finite $Cov_P(IF(\xi))$. By the law of large numbers and Slutsky's theorem, we have $\hat x_n^\lambda\stackrel{p}{\to}x^*$ as $n\to\infty$ and $\lambda\to0$ (at any rate). Therefore $r(\hat x_n^\lambda)\stackrel{p}{\to}r(x^*)$. Moreover, $\sqrt n\langle IF(\xi),\hat P_n-P\rangle=\sqrt n\left(\frac{1}{n}\sum_{i=1}^nIF(\xi_i)-E_P[IF(\xi)]\right)\Rightarrow Y$ where $Y$ is a Gaussian vector with mean 0 and covariance $Cov_P(IF(\xi))$ by the CLT. Thus 
\begin{equation}
\hat x_n^\lambda-x^*=O_p\left(\frac{1}{\sqrt n}\right)+\lambda K+o_p\left(\frac{1}{\sqrt n}+\lambda\right)=O_p\left(\frac{1}{\sqrt n}+\lambda\right)\label{interim4}
\end{equation}
and hence
$$\|\hat x_n^\lambda-x^*\|^2=O_p\left(\frac{1}{n}+\lambda^2\right)$$
which gives
$$r(\hat x_n^\lambda)\|\hat x_n^\lambda-x^*\|^2=o_p\left(\frac{1}{n}+\lambda^2\right)$$
Thus, \eqref{interim2} becomes
\begin{equation}
Z(\hat x_n^\lambda)-Z(x^*)=\frac{1}{2}(\hat x_n^\lambda-x^*)^\top\nabla^2Z(x^*)(\hat x_n^\lambda-x^*)+o_p\left(\frac{1}{n}+\lambda^2\right)\label{interim3}
\end{equation}
Now putting in \eqref{asymptotic expansion}, we get
\begin{eqnarray}
&&Z(\hat x_n^\lambda)-Z(x^*)\notag\\
&=&\frac{1}{2}\left(\langle IF(\xi),\hat P_n-P\rangle+\lambda K+o_p\left(\frac{1}{\sqrt n}+\lambda\right)\right)^\top\nabla^2Z(x^*)\left(\langle IF(\xi),\hat P_n-P\rangle+\lambda K+o_p\left(\frac{1}{\sqrt n}+\lambda\right)\right){}\notag\\
&&{}+o_p\left(\frac{1}{n}+\lambda^2\right)\label{interim5}
\end{eqnarray}
which can be further written as
\begin{equation*}
\frac{1}{2}(\langle IF(\xi),\hat P_n-P\rangle+\lambda K)^\top\nabla^2Z(x^*)(\langle IF(\xi),\hat P_n-P\rangle+\lambda K)+o_p\left(\frac{1}{n}+\lambda^2\right)
\end{equation*}
or
\begin{equation*}
\frac{1}{2}\langle IF(\xi),\hat P_n-P\rangle^\top\nabla^2Z(x^*)\langle IF(\xi),\hat P_n-P\rangle+\frac{1}{2}\lambda^2K^\top\nabla^2Z(x^*)K+\lambda K^\top\nabla^2Z(x^*)\langle IF(\xi),\hat P_n-P\rangle+o_p\left(\frac{1}{n}+\lambda^2\right)
\end{equation*}
\hfill\Halmos\endproof

\proof{Proof of Proposition \ref{trichotomy}.}
By Assumption \ref{solution expansion} and the definition of $\langle IF(\xi),\hat P_n-P\rangle$ in \eqref{IF}, the CLT implies that
$$\sqrt n\langle IF(\xi),\hat P_n-P\rangle\Rightarrow Y$$
where $Y$ is a Gaussian vector with mean 0 and covariance $Cov_P(IF(\xi))$. From \eqref{prelim expression}, and using the continuity of the quadratic function and Slutsky's theorem, we have, when $a=0$ so that $\sqrt n\lambda_n\to0$, that
\begin{eqnarray*}
n\mathcal G(\hat x_n^{\lambda_n})&=&\frac{1}{2}\left(\sqrt n\langle IF(\xi),\hat P_n-P\rangle\right)^\top\nabla^2Z(x^*)\left(\sqrt n\langle IF(\xi),\hat P_n-P\rangle\right)+\frac{1}{2}n\lambda_n^2K^\top\nabla^2Z(x^*)K{}\\
&&{}+\sqrt n\lambda_nK^\top\nabla^2Z(x^*)\left(\sqrt n\langle IF(\xi),\hat P_n-P\rangle\right)+o_p\left(1+n\lambda_n^2\right)\\
&\Rightarrow&\frac{1}{2}Y^\top\nabla^2Z(x^*)Y
\end{eqnarray*}

Likewise, when $a=\infty$ or $-\infty$, which means $\sqrt n\lambda_n\to\infty$ or $-\infty$, we have
\begin{eqnarray*}
\frac{1}{\lambda^2}\mathcal G(\hat x_n^{\lambda_n})&=&\frac{1}{2n\lambda_n^2}\left(\sqrt n\langle IF(\xi),\hat P_n-P\rangle\right)^\top\nabla^2Z(x^*)\left(\sqrt n\langle IF(\xi),\hat P_n-P\rangle\right)+\frac{1}{2}K^\top\nabla^2Z(x^*)K{}\\
&&+\frac{1}{\sqrt n\lambda_n}K^\top\nabla^2Z(x^*)\left(\sqrt n\langle IF(\xi),\hat P_n-P\rangle\right)+o_p\left(\frac{1}{n\lambda_n^2}+1\right)\\
&\Rightarrow&\frac{1}{2}K^\top\nabla^2(x^*)K
\end{eqnarray*}

Moreover, when $\sqrt n\lambda_n\to a$ for some $0<|a|<\infty$, we have 
\begin{eqnarray*}
n\mathcal G(\hat x_n^{\lambda_n})&=&\frac{1}{2}\left(\sqrt n\langle IF(\xi),\hat P_n-P\rangle\right)^\top\nabla^2Z(x^*)\left(\sqrt n\langle IF(\xi),\hat P_n-P\rangle\right)+\frac{1}{2}n\lambda_n^2K^\top\nabla^2Z(x^*)K{}\\
&&{}+\sqrt n\lambda_nK^\top\nabla^2Z(x^*)\left(\sqrt n\langle IF(\xi),\hat P_n-P\rangle\right)+o_p\left(1+n\lambda_n^2\right)\\
&\Rightarrow&\frac{1}{2}Y^\top\nabla^2Z(x^*)Y+\frac{1}{2}a^2K^\top\nabla^2Z(x^*)K+aK^\top\nabla^2Z(x^*)Y
\end{eqnarray*}

Finally, when $\lambda=0$, we either use the same line of arguments as above or observe that we are in the case $a=0$, which reduces to 
$$n\mathcal G(\hat x_n^{EO})\Rightarrow\frac{1}{2}Y'\nabla^2Z(x^*)Y$$
This concludes the proposition.\hfill\Halmos
\endproof

\proof{Proof of Proposition \ref{comparisons}.}
When $a=0$, from Proposition \ref{trichotomy} it is trivial to see that the weak limits of $n\mathcal G(\hat x_n^{\lambda_n})$ and $n\mathcal G(\hat x_n^{EO})$ coincide. When $a=\infty$ or $-\infty$, we have $n\mathcal G(\hat x_n^{\lambda_n})=(n\lambda_n^2)(1/\lambda_n^2)\mathcal G(\hat x_n^{\lambda_n})\stackrel{p}{\to}\infty$ since $(1/2)K^\top\nabla^2Z(x^*)K>0$ by Assumption \ref{nondegenerate}. Finally, when $0<|a|<\infty$, we compare the weak limits of $n\mathcal G(\hat x_n^{\lambda_n})$ and $n\mathcal G(\hat x_n^{EO})$ in Proposition \ref{trichotomy}. Note that $\frac{1}{2}a^2K^\top\nabla^2(x^*)K$ is deterministic and positive. Moreover, $aK^\top\nabla^2(x^*)Y$ satisfies
\begin{equation}
E\left[aK^\top\nabla^2Z(x^*)Y\Bigg|\frac{1}{2}Y^\top\nabla^2Z(x^*)Y\right]=0\label{interim10}
\end{equation}
This is because $Y$, as a mean-zero Gaussian vector, is symmetric, i.e., $Y\stackrel{d}{=}-Y$, and thus for any $z\in\mathbb R$, $Y$ given $(1/2)Y^\top\nabla^2Z(x^*)Y=z$ has the same distribution as $-Y$ given $(1/2)(-Y)^\top\nabla^2Z(x^*)(-Y)=z$ or equivalently $(1/2)Y^\top\nabla^2Z(x^*)Y=z$. This implies $E[Y|(1/2)Y^\top\nabla^2Z(x^*)Y]=E[-Y|(1/2)Y^\top\nabla^2Z(x^*)Y]$, which in turn implies $E[Y|(1/2)Y^\top\nabla^2Z(x^*)Y]=0$ and thus \eqref{interim10}. Hence, by Proposition \ref{mean-preserving spread}, we get that $(1/2)Y^\top\nabla^2Z(x^*)Y$ is second-order stochastically dominated by $(1/2)Y^\top\nabla^2Z(x^*)Y+(1/2)a^2K^\top\nabla^2Z(x^*)K+aK^\top\nabla^2(x^*)Y$.\hfill\Halmos
\endproof

\proof{Proof of Theorem \ref{main thm}.}
This follows immediately from Proposition \ref{comparisons}, allowing for the weak limit of $\infty$ and using the subsequence argument in Definition \ref{asymptotic sosd} when needed.\hfill\Halmos
\endproof

\proof{Proof of Theorem \ref{main thm con}.}
Note that Assumption \ref{assumption:Lagrangian} gives 
\begin{equation}
\nabla Z(x^*)=-\sum_{j\in B}\alpha_j^*\nabla g_j(x^*)\label{interim con}
\end{equation}
Moreover, since Assumption \ref{solution expansion} implies $\hat x_n^\lambda-x^*\stackrel{p}{\to}0$ as $n\to\infty$ and $\lambda\to0$, we have, from Assumption \ref{assumption:preservance}, with probability converging to 1,
$$0=g_j(\hat x_n^\lambda)-g_j(x^*)=\nabla g_j(x^*)^\top(\hat x_n^\lambda-x^*)+\frac{1}{2}(\hat x_n^\lambda-x^*)^\top\nabla^2g_j(x^*)(\hat x_n^\lambda-x^*)+o_p(\|\hat x_n^\lambda-x^*\|^2)$$
or
\begin{equation}
-\nabla g_j(x^*)^\top(\hat x_n^\lambda-x^*)=\frac{1}{2}(\hat x_n^\lambda-x^*)^\top\nabla^2g_j(x^*)(\hat x_n^\lambda-x^*)+o_p(\|\hat x_n^\lambda-x^*\|^2)\label{interim con1}
\end{equation}

So, using \eqref{interim con} and \eqref{interim con1}, we get
\begin{eqnarray*}
&&Z(\hat x_n^\lambda)-Z(x^*)\\
&=&\nabla Z(x^*)^\top(\hat x_n^\lambda-x^*)+\frac{1}{2}(\hat x_n^\lambda-x^*)^\top\nabla^2Z(x^*)(\hat x_n^\lambda-x^*)+o_p(\|\hat x_n^\lambda-x^*\|^2)\\
&=&-\sum_{j\in B}\alpha_j^*\nabla g_j(x^*)^\top(\hat x_n^\lambda-x^*)+\frac{1}{2}(\hat x_n^\lambda-x^*)^\top\nabla^2Z(x^*)(\hat x_n^\lambda-x^*)+o_p(\|\hat x_n^\lambda-x^*\|^2)\\
&=&\frac{1}{2}\sum_{j\in B}\alpha_j^*(\hat x_n^\lambda-x^*)^\top\nabla^2g_j(x^*)(\hat x_n^\lambda-x^*)+\frac{1}{2}(\hat x_n^\lambda-x^*)^\top\nabla^2Z(x^*)(\hat x_n^\lambda-x^*)+o_p(\|\hat x_n^\lambda-x^*\|^2)\\
&=&\frac{1}{2}(\hat x_n^\lambda-x^*)^\top\left(\nabla^2Z(x^*)+\sum_{j\in B}\alpha_j^*\nabla^2g_j(x^*)\right)(\hat x_n^\lambda-x^*)+o_p(\|\hat x_n^\lambda-x^*\|^2)
\end{eqnarray*}
It is now clear that the optimality gap $Z(\hat x_n^\lambda)-Z(x^*)$ behaves the same as \eqref{new opt gap CR} except that $\nabla^2Z(x^*)$ is replaced by $\nabla^2Z(x^*)+\sum_{j\in B}\alpha_j^*\nabla^2g_j(x^*)$. Thus the proofs of all theorems go through in the same way as previously except the said change.\hfill\Halmos
\endproof

\section{Conditions and Proof of Theorem \ref{regularization thm}}\label{sec:proof regularization}

We first impose the following regularity conditions on $h$ and $R$:

\begin{assumption}[Regularity conditions on the loss function]
For any $\xi\in\Xi$, 
\begin{enumerate}
\item $h(\cdot,\xi)$ is a.s. twice continuously differentiable in $x\in\mathcal X$.
\item $h(\cdot,\xi)$ is $\mathcal L^1$-Lipschitz, i.e.,
$$|h(x_1,\xi)-h(x_2,\xi)|\leq L_1(\xi)\|x_1-x_2\|\text{\ \ a.s.}$$
for some non-negative function $L_1(\cdot)$ where $E_P[L_1(\xi)]<\infty$.
\item $\nabla h(\cdot,\xi)$ is $\mathcal L^2$-Lipschitz, i.e.,
$$\|\nabla h(x_1,\xi)-\nabla h(x_2,\xi)\|\leq L_2(\xi)\|x_1-x_2\|\text{\ \ a.s.}$$
for some non-negative function $L_2(\cdot)$ where $E_P[L_2(\xi)^2]<\infty$.
\item $E_P[\nabla^2h(x^*,\xi)]$ is finite and positive definite.
\end{enumerate}\label{regularity loss}
\end{assumption}

\begin{assumption}[Loss function complexity]
$\{\nabla h(x,\cdot):\|x-x^*\|\leq\delta\}$ is Donsker for some small $\delta>0$, and $\{\nabla h(x,\cdot):x\in\mathcal X\}$ is Glivenko-Cantelli.\label{Donsker}
\end{assumption}

\begin{assumption}[Regularity conditions on the penalty function]
We have:
\begin{enumerate}
\item $R$ is twice continuously differentiable in $x\in\mathcal X$.
\item $\nabla^2R(x)$ is uniformly bounded for any $x\in\mathcal X$.
\end{enumerate}\label{regularity penalty}
\end{assumption}

Assumption \ref{regularity loss} stipulates that $h$ is sufficiently smooth, with gradients satisfying Lipschitzness and curvature properties. Assumption \ref{Donsker} is a condition on the ``functional complexity" of the class $\{h(x,\cdot)\}$ and implies that $h(x,\xi)$ satisfies a suitable uniform law of large numbers and CLT, properties that are commonly used to ensure the validity of EO. Assumption \ref{regularity penalty} stipulates that $R$ is sufficiently smooth with a bounded Hessian matrix.

Next we state explicitly that the oracle solution $x^*$ and data-driven solution $\hat x_n^{Reg,\lambda}$ both satisfy corresponding first-order optimality conditions.
\begin{assumption}[First-order optimality conditions]
We have:
\begin{enumerate}
\item $x^*$ is a solution to $\nabla E_P[h(x,\xi)]=0$, and
$\inf_{x\in\mathcal X:\|x-x^*\|\geq\epsilon}\|\nabla E_P[h(x,\xi)]\|>0$ for every $\epsilon>0$.
\item $\hat x_n^{Reg,\lambda}$ is a solution to $\nabla E_{\hat P_n}[h(x,\xi)]+\lambda\nabla R(x)=0$.
\end{enumerate}
\label{first-order assumption}
\end{assumption}
Assumption \ref{first-order assumption} part 1 signifies that the oracle solution $x^*$ satisfies the first-order optimality condition, and moreover is unique in this satisfaction. Part 2 of the assumption signifies that the regularized data-driven solution $\hat x_n^{Reg,\lambda}$ satisfies the first-order opitmality condition of the regularized EO problem.

We are now ready to prove Theorem \ref{regularization thm}. 
\proof{Proof of Theorem \ref{regularization thm}.}
We verify all the conditions in Theorem \ref{nuisance} to conclude the result, where we set $x_n$ in Theorem \ref{nuisance} to be $\hat x_n^{Reg,\lambda_n}$. First of all, by Assumptions \ref{regularity loss}.1 and \ref{regularity loss}.2, we can exchange derivative and expectation so that $\nabla E_P[h(x,\xi)]=E_P[\nabla h(x,\xi)]$ (e.g., Chapter VII, Proposition 2.3 in \citealt{asmussen2007stochastic}). Now, we take $\psi_{x,\lambda}(\xi)$ in Theorem \ref{nuisance} as $\nabla h(x,\xi)+\lambda\nabla R(x)$, and set $\lambda_0=0$ therein. Then:
\\

\noindent\emph{Condition 1: }Assumption \ref{Donsker} together with that $R(x)$ is deterministic and satisfies Assumption \ref{regularity penalty}.1 give that $\{\psi_{x,\lambda}(\cdot):\|x-x^*\|\leq\delta,|\lambda-\lambda_0|\leq\delta\}$ is Donsker (Example 2.10.7 in \citealt{van1996weak}). Thus Condition 1 in Theorem \ref{nuisance} is satisfied.
\\

\noindent\emph{Condition 2: }Consider
\begin{align*}
E_P\|\nabla h(x,\xi)+\lambda\nabla R(x)-\nabla h(x^*,\xi)\|^2&\leq((E_P\|\nabla h(x,\xi)-\nabla h(x^*,\xi)\|^2)^{1/2}+\lambda E_P\|\nabla R(x)\|)^2\\
&\leq((E_P[L_2(\xi)^2])^{1/2}\|x-x^*\|+\lambda E_P\|\nabla R(x)\|)^2\\
&\to0
\end{align*}
as $(x,\lambda)\to(x^*,0)$, where we have used Minkowski's inequality in the first inequality and Assumption \ref{regularity loss}.3 in the second inequality. Thus Condition 2 in Theorem \ref{nuisance} is satisfied.
\\

\noindent\emph{Condition 3: }Condition 3 in Theorem \ref{nuisance} is satisfied since we have argued for Condition 1 above that $\nabla E_P[h(x,\xi)]=E_P[\nabla h(x,\xi)]$, and that $x^*$ satisfies the first-order condition in Assumption \ref{first-order assumption}.1.
\\

\noindent\emph{Condition 4: }Assumptions \ref{regularity loss}.1 and \ref{regularity loss}.3 give that
$$\left\|\frac{\nabla h(x^*+\delta e_j,\xi)-\nabla h(x^*,\xi)}{\delta}\right\|\leq L_2(\xi)$$
 in a neighborhood of 0, where $E_P[L_2(\xi)]<\infty$ and $e_j$ denotes a vector of 0's except the $j$-th component being 1, and $\delta\neq0$. Thus the dominated convergence theorem gives $\nabla E_P[\nabla h(x^*,\xi)]=E_P[\nabla^2h(x^*,\xi)]$. This shows that $E_P[\nabla h(x,\xi)]$ is differentiable at  $x=x^*$. Moreover, together with Assumption \ref{regularity penalty}.1 we have
$$\nabla(E_P[\nabla h(x^*,\xi)]+\lambda\nabla R(x^*))=E_P[\nabla^2 h(x^*,\xi)]+\lambda\nabla^2R(x^*)$$ Furthermore, clearly
$$E_P[\nabla^2 h(x^*,\xi)]+\lambda\nabla^2R(x^*)\to E_P[\nabla^2h(x^*,\xi)]$$
as $\lambda\to0$, and with Assumption \ref{regularity loss}.4 we have $E_P[\nabla^2h(x^*,\xi)]+\lambda\nabla R(x^*)$ being positive definite for $\lambda$ in a neighborhood of 0. Thus Condition 4 in Theorem \ref{nuisance} is satisfied.
\\

\noindent\emph{Condition 5: }Condition 5 in Theorem \ref{nuisance} is satisfied since $\nabla(E_{\hat P_n}[h(x,\xi)]+\lambda R(x))=E_{\hat P_n}[\nabla h(x,\xi)]+\lambda\nabla R(x)$ by Assumptions \ref{regularity loss}.1 and \ref{regularity penalty}.1, and by Assumption \ref{first-order assumption}.2 we have $\hat x_n^{Reg,\lambda}$ being a solution to $E_{\hat P_n}[\nabla h(x,\xi)+\lambda\nabla R(x)]=0$.
\\

\noindent\emph{Condition 6: }We show that Condition 6 in Theorem \ref{nuisance} holds, by setting $x_n$ to be $\hat x_n^{Reg,\lambda_n}$ in Theorem \ref{consistency} and proving its convergence in probability to $x^*$. To this end, consider $\Psi_n(x)=E_{\hat P_n}[\nabla h(x,\xi)]+\lambda_n\nabla R(x)$ and $\Psi(x)=E_P[\nabla h(x,\xi)]$, so that we have $\Psi_n(\hat x_n^{Reg,\lambda_n})=0$ and $\Psi(x^*)=0$ as argued before. Since $\{\nabla h(x,\cdot):x\in\mathcal X\}$ is Glivenko-Cantelli by Assumption \ref{Donsker}, we have
$$\sup_{x\in\mathcal X}\|E_{\hat P_n}[\nabla h(x,\xi)]+\lambda_n\nabla R(x)-E_P[\nabla h(x,\xi)]\|\leq\sup_{x\in\mathcal X}\|E_{\hat P_n}[\nabla h(x,\xi)]-E_P[\nabla h(x,\xi)]\|+\lambda_n\sup_{x\in\mathcal X}\|\nabla R(x)\|\stackrel{p}{\to}0$$
as $n\to\infty$ by Minkowski's equality and Assumption \ref{regularity penalty}.2. This concludes $\sup_{x\in\mathcal X}\|\Psi_n(x)-\Psi(x)\|\stackrel{p}{\to}0$ in Theorem \ref{consistency}. Moreover, Assumption \ref{first-order assumption} concludes 
$\inf_{x\in\mathcal X:\|x-x^*\|\geq\epsilon}\|\Psi(x)\|>0=\|\Psi(x^*)\|$. Therefore, by Theorem \ref{consistency}, we have $\hat x_n^{Reg,\lambda_n}\stackrel{p}{\to}x^*$.
\\

With the above six conditions, we can invoke Theorem \ref{nuisance} to obtain
\begin{eqnarray*}
&&\hat x_n^{Reg,\lambda_n}-x^*\\
&=&-(E_P[\nabla^2h(x^*,\xi)])^{-1}E_P[\nabla h(x^*,\xi)+\lambda_n\nabla R(x)]-(E_P[\nabla^2h(x^*,\xi)])^{-1} E_{\hat P_n}[\nabla h(x^*,\xi)]{}\\
&&{}+o_p\left(\frac{1}{\sqrt n}+\|E[\nabla h(x^*,\xi)+\lambda_n\nabla R(x^*)]\|\right)\\
&=&-\langle (E_P[\nabla^2h(x^*,\xi)])^{-1}\nabla h(x^*,\xi),\hat P_n-P\rangle-\lambda_n(E_P[\nabla^2h(x^*,\xi)])^{-1}\nabla R(x^*)+o_p\left(\frac{1}{\sqrt n}+\lambda_n\|\nabla R(x^*)\|\right)
\end{eqnarray*}
since $E_P[\nabla h(x^*,\xi)]=0$ by the first-order optimality condition in Assumption \ref{first-order assumption}.1.\hfill\Halmos
\endproof

\section{Conditions and Proofs of Theorem \ref{thm:main DRO}}\label{sec:DRO proof}

We need the following technical assumption on the involved optimization problem.

\begin{assumption}[Validity of KKT conditions]
With probability converging to 1 as $n\to\infty$, we have
\begin{enumerate}
\item The optimization problem $\min_{x\in\mathcal X,\alpha\geq0,\beta\in\mathbb R}\left\{\alpha E_{\hat P_n}\left[\phi^*\left(\frac{h(x,\xi)-\beta}{\alpha}\right)\right]+\alpha\lambda_n+\beta\right\}$ has a unique solution $(\hat x_n^*,\hat\alpha_n^*,\hat\beta_n^*)$ that satisfies the KKT conditions with $\hat x_n^*\in\mathcal X^\circ$, the interior of $\mathcal X$.
\item For all $x\in\mathcal X$, $\max_{D(Q,\hat P_n)\leq\lambda_n}E_Q[h(x,\xi)]\neq\text{ess sup}_{\hat P_n}\ h(x,\xi)$.
\end{enumerate}\label{KKT}
\end{assumption}

The first condition in Assumption \ref{KKT} guarantees that we can use the KKT conditions to locate the optimal solution of \eqref{DRO} (we will soon see that the problem $\min_{x\in\mathcal X,\alpha\geq0,\beta\in\mathbb R}\left\{\alpha E_{\hat P_n}\left[\phi^*\left(\frac{h(X)-\beta}{\alpha}\right)\right]+\alpha\lambda_n+\beta\right\}$ is a dualization of \eqref{DRO}). The second condition ensures the obtained optimal solution is nontrivial. Both conditions need to be satisfied with overwhelming probability as $n\to\infty$.

Analogous to Assumption \ref{first-order assumption}, we also need the following first-order optimality condition for the oracle solution.
\begin{assumption}[First-order optimality conditions]
$x^*$ is a solution to $\nabla E_P[h(x,\xi)]=0$, and
$\inf_{x\in\mathcal X:\|x-x^*\|\geq\epsilon}\|\nabla E_P[h(x,\xi)]\|>0$ for every $\epsilon>0$.
\label{first-order assumption DRO}
\end{assumption}


Next we present assumptions on $h$.

\begin{assumption}[Regularity conditions on the loss function]
For any $\xi\in\Xi$, 
\begin{enumerate}
\item $h(\cdot,\xi)$ is a.s. twice continuously differentiable in $x\in\mathcal X$.
\item $h(\cdot,\cdot)$ and $\nabla h(\cdot,\cdot)$ are uniformly bounded over $x\in\mathcal X$ and $\xi\in\Xi$.
\item $\nabla h(\cdot,\xi)$ is $\mathcal L^2$-Lipschitz, i.e.,
$$\|\nabla h(x_1,\xi)-\nabla h(x_2,\xi)\|\leq L_2(\xi)\|x_1-x_2\|\text{\ \ a.s.}$$
for some non-negative function $L_2(\cdot)$ where $E_P[L_2(\xi)^2]<\infty$. 
\item $E_P[\nabla^2h(x^*,\xi)]$ is finite and positive definite.
\end{enumerate}\label{regularity DRO}
\end{assumption}


Compared to Assumption \ref{regularity loss}, here we strengthen part 2 of the assumption to the uniform boundedness of $h$ and $\nabla h$. Note that the uniform boundedness of $\nabla h$ implies $\mathcal L^1$-Lipschitzness of $h$, which is precisely Assumption \ref{regularity loss}.2. The stronger Assumption \ref{regularity DRO}.2 can be potentially relaxed but it helps streamline our mathematical arguments.

\begin{assumption}[Loss function complexity]
We have:
\begin{enumerate}
\item $\{h(x,\cdot):x\in\mathcal X\}$, $\{h(x,\cdot)^2:x\in\mathcal X\}$ and $\{\nabla h(x,\cdot):x\in\mathcal X\}$ are Glivenko-Cantelli.
\item $\{{\phi^*}'(\alpha h(x,\cdot)-\beta)\nabla h(x,\cdot):\|x-x^*\|\leq\delta,0\leq\alpha\leq\delta,0\leq\beta\leq\delta\}$ is Donsker for some small $\delta>0$.
\item $\inf_{x\in\mathcal X}Var_P(h(x,\xi))>0$.
\end{enumerate}
\label{Donsker DRO}
\end{assumption}

Compared to Assumption \ref{Donsker}, Assumption \ref{Donsker DRO} requires the Glivenko-Cantelli property for $h$ and $h^2$, and the Donsker property for a larger function class $\{{\phi^*}'(\alpha h(x\cdot)-\beta)\nabla h(x,\cdot)\}$. Moreover, it requires a non-degeneracy on $h$ uniformly across all $x\in\mathcal X$.

We are now ready to prove Theorem \ref{thm:main DRO}.
\proof{Proof of Theorem \ref{thm:main DRO}.}
Consider the inner maximization in \eqref{DRO}. By a change of measure from $Q$ to $\hat P_n$, we rewrite it as
\begin{equation}
\begin{array}{ll}
\max_{L\geq0}&E_{\hat P_n}[h(x,\xi)L]\\
\text{subject to}&E_{\hat P_n}[\phi(L)]\leq\lambda
\end{array}\label{max DRO transformed}
\end{equation}
where the decision variable is now the likelihood ratio $L=L(\xi)$. Next, since $L\equiv1$ is in the interior of the feasible region for any $\lambda>0$, Slater's condition holds and we rewrite \eqref{max DRO transformed} as the Lagrangian formulation
\begin{eqnarray}
&&\min_{\alpha\geq0,\beta}\max_{L\geq0}E_{\hat P_n}[h(x,\xi)L]-\alpha(E_{\hat P_n}[\phi(L)]-\lambda)-\beta(E_{\hat P_n}[L]-1)\notag\\
&=&\min_{\alpha\geq0,\beta}\max_{L\geq0}E_{\hat P_n}[(h(x,\xi)-\beta)L-\alpha\phi(L)]+\alpha\lambda+\beta\notag\\
&=&\min_{\alpha\geq0,\beta}E_{\hat P_n}\left[\max_{L\geq0}\{(h(x,\xi)-\beta)L-\alpha\phi(L)\}\right]+\alpha\lambda+\beta\notag\\
&=&\min_{\alpha\geq0,\beta}\alpha E_{\hat P_n}\left[\max_{L\geq0}\left\{\frac{h(x,\xi)-\beta}{\alpha}L-\phi(L)\right\}\right]+\alpha\lambda+\beta\notag\\
&=&\min_{\alpha\geq0,\beta}\alpha E_{\hat P_n}\left[\phi^*\left(\frac{h(x,\xi)-\beta}{\alpha}\right)\right]+\alpha\lambda+\beta\label{interim DRO}
\end{eqnarray}
where in the above we define $0\phi^*(a/0)=\infty$ if $a>0$ and $0\phi^*(a/0)=0$ for $a\leq0$ (i.e., when $\alpha=0$ the resulting objective function is equal to $\infty$ if $\text{ess sup}_{\hat P}\ h(x,\xi)>\beta$, and $\beta$ otherwise). 
Thus, \eqref{DRO} is equivalent to
\begin{equation}
\min_{x\in\mathcal X,\alpha\geq0,\beta}\alpha E_{\hat P_n}\left[\phi^*\left(\frac{h(x,\xi)-\beta}{\alpha}\right)\right]+\alpha\lambda+\beta\label{opt dual new}
\end{equation}

Now suppose the two events in Assumption \ref{KKT} hold, which occur with probability reaching 1 as $n\to\infty$. In the first event, the KKT conditions of \eqref{opt dual new} uniquely identifies the optimal solution. Denote $(\hat x_n^\lambda,\hat\alpha_n^\lambda,\hat\beta_n^\lambda)$ as the optimal solution of \eqref{opt dual new}. If $\hat\alpha_n^\lambda=0$, then $\hat\beta_n^\lambda=\text{ess sup}_{\hat P_n}\ h(\hat x_n^\lambda,\xi)$ and the objective value of \eqref{opt dual new} becomes $\text{ess sup}_{\hat P_n}\ h(\hat x_n^\lambda,\xi)$, which violates the second event. Thus, with probability reaching 1, we must have $\hat\alpha_n^\lambda>0$. We will focus on this case for the rest of the proof and use a standard convergence in probability argument to handle the other case (which occurs with vanishing probability) when needed.

The KKT conditions in the case $\hat\alpha_n^\lambda>0$ are precisely the first-order optimality conditions
\begin{eqnarray*}
&&E_{\hat P_n}\left[{\phi^*}'\left(\frac{h(x,\xi)-\beta}{\alpha}\right)\nabla h(x,\xi)\right]=0\\
&&E_{\hat P_n}\left[\phi^*\left(\frac{h(x,\xi)-\beta}{\alpha}\right)\right]-E_{\hat P_n}\left[{\phi^*}'\left(\frac{h(x,\xi)-\beta}{\alpha}\right)\frac{h(x,\xi)-\beta}{\alpha}\right]+\lambda=0\\
&&-E_{\hat P_n}\left[{\phi^*}'\left(\frac{h(x,\xi)-\beta}{\alpha}\right)\right]+1=0
\end{eqnarray*}
or equivalently
\begin{eqnarray*}
&&E_{\hat P_n}\left[{\phi^*}'\left(\frac{h(x,\xi)-\beta}{\alpha}\right)\nabla h(x,\xi)\right]=0\\
&&E_{\hat P_n}\left[{\phi^*}'\left(\frac{h(x,\xi)-\beta}{\alpha}\right)\frac{h(x,\xi)-\beta}{\alpha}\right]-E_{\hat P_n}\left[\phi^*\left(\frac{h(x,\xi)-\beta}{\alpha}\right)\right]=\lambda\\
&&E_{\hat P_n}\left[{\phi^*}'\left(\frac{h(x,\xi)-\beta}{\alpha}\right)\right]=1
\end{eqnarray*}

For convenience, we do a change of variable and define $\tilde\alpha=1/\alpha$ and $\tilde\beta=\beta/\alpha$. Then the above can be written as
\begin{eqnarray}
&&E_{\hat P_n}\left[{\phi^*}'(\tilde\alpha h(x,\xi)-\tilde\beta)\nabla h(x,\xi)\right]=0\label{interim DRO new}\\
&&E_{\hat P_n}\left[{\phi^*}'(\tilde\alpha h(x,\xi)-\tilde\beta)(\tilde\alpha h(x,\xi)-\tilde\beta)\right]-E_{\hat P_n}\left[\phi^*(\tilde\alpha h(x,\xi)-\tilde\beta)\right]=\lambda\label{interim DRO1}\\
&&E_{\hat P_n}\left[{\phi^*}'(\tilde\alpha h(x,\xi)-\tilde\beta)\right]=1\label{interim DRO2}
\end{eqnarray}
and we denote $(\hat x_n^\lambda,\tilde\alpha_n^\lambda,\tilde\beta_n^\lambda)$ as the root to \eqref{interim DRO new}-\eqref{interim DRO2}.

We verify Conditions 1-6 in Theorem \ref{nuisance}. The main idea is to view $(\tilde\alpha_n^{\lambda_n},\tilde\beta_n^{\lambda_n})$ as the estimated parameter in Theorem \ref{nuisance} (i.e., the $\lambda_n$ in Theorem \ref{nuisance} is now regarded as $(\tilde\alpha_n^{\lambda_n},\tilde\beta_n^{\lambda_n})$). Then $\hat x_n^{\lambda_n}$ is viewed as the root to  \eqref{interim DRO new}, which is parametrized by $(\tilde\alpha_n^{\lambda_n},\tilde\beta_n^{\lambda_n})$ that in turn depends on $\lambda_n$. In the following, we first show Condition 6 in Theorem \ref{nuisance}, in the form $(\hat x_n^{\lambda_n},\tilde\alpha_n^{\lambda_n},\tilde\beta_n^{\lambda_n})\stackrel{p}{\to}(x^*,0,0)$ as $n\to\infty$ and $\lambda_n\to0$ (at any rate).
\\

\noindent\emph{Condition 6: }We derive the asymptotic behaviors of $\tilde\alpha_n^\lambda$ and $\tilde\beta_n^\lambda$. 
First consider \eqref{interim DRO2}. We show that we can find a unique root $\tilde\beta_n=\tilde\beta_n(\tilde\alpha)$ for \eqref{interim DRO2} as a function of $\tilde\alpha$, at a neighborhood of $\tilde\alpha=0$. Here and from now on, we note that the values of $\tilde\alpha$ and $\tilde\beta$ could depend on $x$ but we will suppress this dependence unless needed. Note that $\tilde\alpha=\tilde\beta=0$ solves \eqref{interim DRO2}. Moreover, by Lemma \ref{phi lemma} $E_{\hat P_n}[{\phi^*}'(\tilde\alpha h(x,\xi)-\tilde\beta)]-1$ is continuously differentiable in $\tilde\alpha$ and $\tilde\beta$ both at 0 and
$$\frac{\partial}{\partial\tilde\beta}E_{\hat P_n}[{\phi^*}'(\tilde\alpha h(x,\xi)-\tilde\beta)]-1\Bigg|_{\tilde\alpha=\tilde\beta=0}=-E_{\hat P_n}[{\phi^*}''(\tilde\alpha h(x,\xi)-\tilde\beta)]\Bigg|_{\tilde\alpha=\tilde\beta=0}=-{\phi^*}''(0)<0$$
So by the implicit function theorem, we have $\tilde\beta_n(\tilde\alpha)$ continuously differentiable in a neighborhood around $\tilde\alpha=0$. Furthermore, since
$$\frac{\partial}{\partial\tilde\alpha}E_{\hat P_n}[{\phi^*}'(\tilde\alpha h(x,\xi)-\tilde\beta)]-1\Bigg|_{\tilde\alpha=\tilde\beta=0}=E_{\hat P_n}[{\phi^*}''(\tilde\alpha h(x,\xi)-\tilde\beta)h(x,\xi)]\Bigg|_{\tilde\alpha=\tilde\beta=0}={\phi^*}''(0)\hat E[h(x,\xi)]$$
we have $\tilde\beta'(0)=-(-{\phi^*}''(0)^{-1}){\phi^*}''(0)E_{\hat P_n}[h(x,\xi)]=\hat E[h(x,\xi)]$. Moreover, since $h$ is uniformly bounded by Assumption \ref{regularity DRO}.2, we have $\tilde\alpha h(x,\xi)-\tilde\beta\to0$ uniformly over $x,\xi$ when $\tilde\alpha,\tilde\beta\to0$. Thus, together with the mean value theorem, we have 
\begin{equation}
\tilde\beta_n(\tilde\alpha)=\tilde\beta_n'(\zeta)\tilde\alpha\label{interim DRO proof3}
\end{equation}
for some $\zeta$ between 0 and $\tilde\alpha$, when $\tilde\alpha$ is in a sufficiently small neighborhood of 0, where $\tilde\beta_n'(\zeta)$ converges to $E_{\hat P_n}[h(x,\xi)]$ uniformly over $\hat P_n$ and $x\in\mathcal X$ as $\tilde\alpha\to0$.

Now consider \eqref{interim DRO1}. By the mean value theorem, we can write it as
\begin{eqnarray}
&&E_{\hat P_n}[{\phi^*}'(0)(\tilde\alpha h(x,\xi)-\tilde\beta)]+E_{\hat P_n}[{\phi^*}''(\zeta_1(x,\xi))(\tilde\alpha h(x,\xi)-\tilde\beta)^2]-E_{\hat P_n}[{\phi^*}'(0)(\tilde\alpha h(x,\xi)-\tilde\beta)]{}\notag\\
&&{}-\frac{1}{2}E_{\hat P_n}[{\phi^*}''(\zeta_2(x,\xi))(\tilde\alpha h(x,\xi)-\tilde\beta)^2]=\lambda\label{interim DRO3}
\end{eqnarray}
which gives
\begin{equation}
E_{\hat P_n}[{\phi^*}''(\zeta_1(x,\xi))(\tilde\alpha h(x,\xi)-\tilde\beta)^2]-\frac{1}{2}E_{\hat P_n}[{\phi^*}''(\zeta_2(x,\xi))(\tilde\alpha h(x,\xi)-\tilde\beta)^2]=\lambda\label{interim DRO3}
\end{equation}
where $\zeta_1(x,\xi)$ and $\zeta_2(x,\xi)$ are some values between $0$ and $\tilde\alpha h(x,\xi)-\tilde\beta$. Let us take $\tilde\beta=\tilde\beta_n(\tilde\alpha)$ as the root of \eqref{interim DRO2} described above, which is well-defined for $\tilde\alpha$ in a neighborhood of 0. We can write the left hand side of \eqref{interim DRO3} as
$$E_{\hat P_n}[{\phi^*}''(\zeta_1(x,\xi))(\tilde\alpha h(x,\xi)-\tilde\beta_n(\tilde\alpha))^2]-\frac{1}{2}E_{\hat P_n}[{\phi^*}''(\zeta_2(x,\xi))(\tilde\alpha h(x,\xi)-\tilde\beta_n(\tilde\alpha))^2]$$
which can be written further as
$$\left(E_{\hat P_n}[{\phi^*}''(\zeta_1(x,\xi))(h(x,\xi)-\tilde\beta_n'(\zeta))^2]-\frac{1}{2}E_{\hat P_n}[{\phi^*}''(\zeta_2(x,\xi))(h(x,\xi)-\tilde\beta_n'(\zeta))^2]\right)\tilde\alpha^2$$
where $\zeta$ lies between 0 and $\tilde\alpha$. For further convenience, we denote
$$u_n(\tilde\alpha)=E_{\hat P_n}[{\phi^*}''(\zeta_1(x,\xi))(h(x,\xi)-\tilde\beta_n'(\zeta))^2]-\frac{1}{2}E_{\hat P_n}[{\phi^*}''(\zeta_2(x,\xi))(h(x,\xi)-\tilde\beta_n'(\zeta))^2]$$
Note that $u_n(\tilde\alpha)$ is continuous at $\tilde\alpha=0$, and $u_n(0)={\phi^*}''(0)Var_{\hat P_n}(h(x,\xi))/2$. We argue that $u_n(\tilde\alpha)$ converges to ${\phi^*}''(0)Var_P(h(x,\xi))/2>0$ uniformly over $x\in\mathcal X$ a.s. as $n\to\infty$ and $\tilde\alpha\to0$ (at any rate relative to $n$). More precisely, we consider 
\begin{equation}
\sup_{x\in\mathcal X}\left|u_n(\tilde\alpha)-\frac{{\phi^*}''(0)Var_P(h(x,\xi))}{2}\right|\leq\sup_{x\in\mathcal X}|u_n(\tilde\alpha)-u_n(0)|+\sup_{x\in\mathcal X}\left|u_n(0)-\frac{{\phi^*}''(0)Var_P(h(x,\xi))}{2}\right|\label{interim uniform}
\end{equation}
Now, by using the uniform boundedness of $h$ (Assumption \ref{regularity DRO}.2) and our conclusion above that $\tilde\beta_n'(\zeta)\to E_{\hat P_n}[h(x,\xi)]$ uniformly over $\hat P_n,x$ as $\tilde\alpha\to0$, we obtain $\sup_{x\in\mathcal X}|u_n(\tilde\alpha)-u_n(0)|\to0$. On the other hand, by the Glivenko-Cantelli properties of the classes of $h$ and $h^2$ (Assumption \ref{Donsker DRO}.1), we have
$$\sup_{x\in\mathcal X}\left|u_n(0)-\frac{{\phi^*}''(0)Var_P(h(x,\xi))}{2}\right|=\frac{{\phi^*}''(0)}{2}\sup_{x\in\mathcal X}\left|Var_{\hat P_n}(h(x,\xi))-Var_P(h(x,\xi))\right|\to0\text{\ \ a.s.}$$
Thus \eqref{interim uniform} goes to 0. Hence, given a sufficiently small $\lambda$ and large $n$, we can find a root $\tilde\alpha_n^\lambda$ to \eqref{interim DRO3}, or equivalently
$$u_n(\tilde\alpha_n^\lambda)(\tilde\alpha_n^\lambda)^2=\lambda$$
in a small neighborhood of 0, uniformly over all $x\in\mathcal X$. Moreover, note that the events represented in Assumption \ref{KKT} occurs with probability tending to 1, and thus, by the uniform non-degeneracy of the variance (Assumption \ref{Donsker DRO}.3), we have $\tilde\alpha_n^\lambda$ satisfy
\begin{equation}
\tilde\alpha_n^\lambda=\sqrt{\frac{\lambda}{u_n(\tilde\alpha_n^\lambda)}}\stackrel{p}{\to}0
\label{interim DRO proof2}
\end{equation}
as $n\to\infty$ and $\lambda\to0$.
Correspondingly, $$\tilde\beta_n^\lambda=\tilde\beta_n(\tilde\alpha_n^\lambda)=\tilde\beta_n'(\zeta_n^\lambda)\tilde\alpha_n^\lambda\stackrel{p}{\to}0
$$
where $\zeta_n^\lambda$ lies between 0 and $\tilde\alpha_n^\lambda$. In other words, $(\tilde\alpha_n^\lambda,\tilde\beta_n^\lambda)\stackrel{p}{\to}0$ as $n\to\infty$ and $\lambda\to0$. 

Next we show that $\hat x_n^\lambda\stackrel{p}{\to}x^*$. To this end, we verify the conditions in Theorem \ref{consistency} where $\Phi(x)$ is set to be $\nabla E_P[h(x,\xi)]=E_P[\nabla h(x,\xi)]$ (by the $\mathcal L^1$-Lipschitzness of $h$ implied by Assumption \ref{regularity DRO}.2) and $\Phi_n(x)$ is set to be $E_{\hat P_n}\left[{\phi^*}'\left(\tilde\alpha_n^{\lambda_n} h(x,\xi)-\tilde\beta_n^{\lambda_n}\right)\nabla h(x,\xi)\right]$ in \eqref{interim DRO new}. We have
$$\sup_{x\in\mathcal X}\left\|E_{\hat P_n}\left[{\phi^*}'\left(\tilde\alpha_n^{\lambda_n} h(x,\xi)-\tilde\beta_n^{\lambda_n}\right)\nabla h(x,\xi)\right]-E_P[\nabla h(x,\xi)]\right\|\stackrel{p}{\to}0$$
by the continuity of ${\phi^*}'$ at 0 (Lemma \ref{phi lemma}), the Glivenko-Cantelli property of $\nabla h$ (Assumption \ref{Donsker DRO}.1), the uniform boundedness of $h$ and $\nabla h$ (Assumption \ref{regularity DRO}.2), and our earlier conclusion that $(\tilde\alpha_n^\lambda,\tilde\beta_n^\lambda)\stackrel{p}{\to}0$. Moreover, $\inf_{x\in\mathcal X:\|x-x^*\|\geq\epsilon}\|\nabla Eh(x,\xi)\|>0$ for every $\epsilon>0$ by Assumption \ref{first-order assumption DRO}. Finally, from our earlier conclusion, we have $\hat x_n^\lambda$ being the root to
$$E_{\hat P_n}\left[{\phi^*}'\left(\tilde\alpha_n^\lambda h(x,\xi)-\tilde\beta_n^\lambda\right)\nabla h(x,\xi)\right]=0$$
with probability tending to 1. Thus, all the conditions in Theorem \ref{consistency} are satisfied and we have $\hat x_n^\lambda\stackrel{p}{\to}x^*$. Therefore, we obtain $\left(\hat x_n^\lambda,\tilde\alpha_n^\lambda,\tilde\beta_n^\lambda\right)\stackrel{p}{\to}0$ as $n\to\infty$ and $\lambda\to0$.

Moreover, from \eqref{interim DRO proof3} and \eqref{interim DRO proof2}, and using the uniform boundedness and continuity of $h$ (Assumptions \ref{regularity DRO}.1 and \ref{regularity DRO}.2), we invoke the dominated convergence theorem to get further that
\begin{equation}
\tilde\alpha_n^\lambda=\sqrt{\frac{2\lambda}{{\phi^*}''(0)Var_P(h(x^*,\xi))}}(1+o_p(1))\label{interim DRO proof5}
\end{equation}
and
\begin{equation}
\tilde\beta_n^\lambda=E_P[h(x^*,\xi)]\sqrt{\frac{2\lambda}{{\phi^*}''(0)Var_P(h(x^*,\xi))}}(1+o_p(1))\label{interim DRO proof6}
\end{equation}
\\

We now continue to verify all other conditions in Theorem \ref{nuisance}. As mentioned, we set the $\lambda$ in the theorem as $(\tilde\alpha,\tilde\beta)$. We then set  $\varphi_{x,\tilde\alpha,\tilde\beta}(\xi)={\phi^*}'(\tilde\alpha h(x,\xi)-\tilde\beta)\nabla h(x,\xi)$ which is the function inside the expectation in \eqref{interim DRO new}.
\\

\noindent\emph{Condition 1: }Direct use of Assumption \ref{Donsker DRO}.2.
\\

\noindent\emph{Condition 2: }
We have
\begin{eqnarray}
&&E_P\|{\phi^*}'(\tilde\alpha h(x,\xi)-\tilde\beta)\nabla h(x,\xi)-\nabla h(x^*,\xi)\|^2\notag\\
&\leq&((E_P\|({\phi^*}'(\tilde\alpha h(x,\xi)-\tilde\beta)-1)\nabla h(x,\xi)\|^2)^{1/2}+(E_P\|\nabla h(x,\xi)-\nabla h(x^*,\xi)\|^2)^{1/2})^2\label{interim DRO proof1}
\end{eqnarray}
by Minkowski's inequality. 
For the first term in \eqref{interim DRO proof1}, since $h$ and $\nabla h$ are uniformly bounded (Assumption \ref{regularity DRO}.2) and ${\phi^*}'$ is continuous at 0 (Lemma \ref{phi lemma}), as $\tilde\alpha,\tilde\beta\to0$, we have
$E_P\|({\phi^*}'(\tilde\alpha h(x,\xi)-\tilde\beta)-1)\nabla h(x,\xi)\|^2\to0$ uniformly in $x$. For the second term in \eqref{interim DRO proof1}, we have
$$E_P\|\nabla h(x,\xi)-\nabla h(x^*,\xi)\|^2\leq E_P[L_2(\xi)^2]\|x-x^*\|^2\to0$$
as $x\to x^*$ since $E[L_2(\xi)^2]<\infty$ (Assumption \ref{regularity DRO}.3). Therefore, putting together we have \eqref{interim DRO proof1} go to 0. This shows the satisfaction of Condition 2.
\\



\noindent\emph{Condition 3: }By the a.s. differentiability of $h$ in Assumption \ref{regularity DRO}.1 and the $\mathcal L^1$-Lipschitzness of $h$ implied by Assumption \ref{regularity DRO}.2, we have $\nabla E_P[h(x,\xi)]=E_P[\nabla h(x,\xi)]$. Thus, Assumption \ref{first-order assumption DRO} concludes the satisfaction of Condition 3.
\\

\noindent\emph{Condition 4: }Since $h$ is uniformly bounded (Assumption \ref{regularity DRO}.2) and also $\mathcal L^1$-Lipschitz (implied by Assumption \ref{regularity DRO}.2), and ${\phi^*}'$ is continuously differentiable at 0 (Lemma \ref{phi lemma}) and thus Lipschitz continuous at a neighborhood of 0, we have ${\phi^*}'(\tilde\alpha h(x,\xi)-\tilde\beta)$ being $\mathcal L^1$-Lipschitz and uniformly bounded over $x,\xi$ a.s. for a fixed sufficiently small $(\tilde\alpha,\tilde\beta)$. On the other hand, we have $\nabla h$ being $\mathcal L^1$-Lipschitz (implied by Assumption \ref{regularity DRO}.3) and uniformly bounded (Assumption \ref{regularity DRO}.2). Thus we have ${\phi^*}'(\tilde\alpha h(x,\xi)-\tilde\beta)\nabla h(x,\xi)$ being $\mathcal L^1$-Lipschitz. Moreover, by Assumption \ref{regularity DRO}.1 and Lemma \ref{phi lemma} again it is differentiable a.s. Thus, using the same argument used to verify Condition 4 in the proof of Theorem \ref{regularization thm}, we can exchange derivative and expectation to obtain
\begin{align*}
\nabla E_P[{\phi^*}'(\tilde\alpha h(x,\xi)-\tilde\beta)\nabla h(x,\xi)]&=E_P[\nabla({\phi^*}'(\tilde\alpha h(x,\xi)-\tilde\beta)\nabla h(x,\xi))]\\
&=E_P[{\phi^*}''(\tilde\alpha h(x,\xi)-\tilde\beta)\tilde\alpha\nabla h(x,\xi)\nabla h(x,\xi)']+E_P[{\phi^*}'(\tilde\alpha h(x,\xi)-\tilde\beta)\nabla^2h(x,\xi)]
\end{align*}
for $(\tilde\alpha,\tilde\beta)$ in a neighborhood of 0. In particular, $\nabla E_P[\nabla h(x,\xi)]=E_P[\nabla^2h(x,\xi)]$. 

Now, by the uniform boundedness of $h$ and $\nabla h$ (Assumption \ref{regularity DRO}.1), the continuity of ${\phi^*}'$ and ${\phi^*}''$ at 0 (Lemma \ref{phi lemma}) and the finiteness of $E_P[\nabla^2h(x^*,\xi)]$ (Assumption \ref{regularity DRO}.4), we use the dominated convergence theorem to obtain
$$E_P[{\phi^*}''(\tilde\alpha h(x^*,\xi)-\tilde\beta)\tilde\alpha\nabla h(x^*,\xi)\nabla h(x^*,\xi)']+E_P[{\phi^*}'(\tilde\alpha h(x^*,\xi)-\tilde\beta)\nabla^2h(x^*,\xi)]
\to E_P[\nabla^2h(x^*,\xi)]$$
as $(\tilde\alpha,\tilde\beta)\to(0,0)$.

Finally, by the continuity of ${\phi^*}'$ ${\phi^*}''$ at 0 (Lemma \ref{phi lemma}), the uniform boundedness of $h$ and $\nabla h$ (Assumption \ref{regularity DRO}.1), and Assumption \ref{regularity DRO}.4, we have $E_P[{\phi^*}''(\tilde\alpha h(x^*,\xi)-\tilde\beta)\tilde\alpha\nabla h(x^*,\xi)\nabla h(x^*,\xi)']+E_P[{\phi^*}'(\tilde\alpha h(x^*,\xi)-\tilde\beta)\nabla^2h(x^*,\xi)]$ being positive definite with eigenvalues uniformly bounded and away from 0 for $(\tilde\alpha,\tilde\beta)$ in a neighborhood of $(0,0)$. This concludes the satisfaction of Condition 4.
\\

\noindent\emph{Condition 5: }Since $\hat x_n^\lambda$ satisfies
$E_{\hat P_n}\left[{\phi^*}'\left(\tilde\alpha_n^\lambda h(x,\xi)-\tilde\beta_n^\lambda\right)\nabla h(x,\xi)\right]=0$ with probability tending to 1, Condition 5 is readily satisfied.
\\

We have therefore verified all conditions in Theorem \ref{nuisance}. Hence,
\begin{eqnarray}
&&\hat x_n^{\lambda_n}-x^*\notag\\
&=&-(E_P[\nabla^2h(x,\xi)])^{-1}E_P\left[{\phi^*}'\left(\tilde\alpha_n^{\lambda_n} h(x^*,\xi)-\tilde\beta_n^{\lambda_n}\right)\nabla h(x^*,\xi)\right]-(E_P[\nabla^2h(x,\xi)])^{-1}E_{\hat P_n}[\nabla h(x^*,\xi)]{}\notag\\
&&{}+o_p\left(\frac{1}{\sqrt n}+\left\|E_P\left[{\phi^*}'\left(\tilde\alpha_n^{\lambda_n} h(x^*,\xi)-\tilde\beta_n^{\lambda_n}\right)\nabla h(x^*,\xi)\right]\right\|\right)\notag\\
&=&-(E_P[\nabla^2h(x,\xi)])^{-1}\left(E_P[\nabla h(x^*,\xi)]+E_P\left[{\phi^*}''(\zeta_n)\left(\tilde\alpha_n^{\lambda_n} h(x^*,\xi)-\tilde\beta_n^{\lambda_n}\right)\nabla h(x^*,\xi)\right]\right){}\notag\\
&&{}-(E_P[\nabla^2h(x,\xi)])^{-1}E_{\hat P_n}[\nabla h(x^*,\xi)]+o_p\left(\frac{1}{\sqrt n}+\left\|E_P\left[{\phi^*}''(\zeta_n)\left(\tilde\alpha_n^{\lambda} h(x^*,\xi)-\tilde\beta_n^{\lambda_n}\right)\nabla h(x^*,\xi)\right]\right\|\right)\notag\\
&=&-(E_P[\nabla^2h(x,\xi)])^{-1}E_P\left[{\phi^*}''(\zeta_n)\left(\tilde\alpha_n^{\lambda_n} h(x^*,\xi)-\tilde\beta_n^{\lambda_n}\right)\nabla h(x^*,\xi)\right]{}\notag\\
&&{}-(E_P[\nabla^2h(x,\xi)])^{-1}E_{\hat P_n}[\nabla h(x^*,\xi)]+o_p\left(\frac{1}{\sqrt n}+\left\|E\left[{\phi^*}''(\zeta_n)\left(\tilde\alpha_n^{\lambda_n} h(x^*,\xi)-\tilde\beta_n^{\lambda_n}\right)\nabla h(x^*,\xi)\right]\right\|\right)\label{interim DRO proof4}
\end{eqnarray}
where $\zeta_n$ lies between 0 and $\tilde\alpha_n^{\lambda_n} h(x^*,\xi)-\tilde\beta_n^{\lambda_n}$, and we have used the first-order optimality condition $E_P[\nabla h(x^*,\xi)]=0$ derived previously. Now, substituting \eqref{interim DRO proof5} and \eqref{interim DRO proof6}, and using the continuity of ${\phi^*}''$ at 0 (Lemma \ref{phi lemma}), uniform boundedness of $h$ and $\nabla h$ (Assumption \ref{regularity DRO}.1) and the characterizations in \eqref{interim DRO proof3} and \eqref{interim DRO proof2}, we invoke the bounded convergence theorem to obtain
\begin{eqnarray*}
&&E_P\left[{\phi^*}''(\zeta_n)\left(\tilde\alpha_n^{\lambda_n} h(x^*,\xi)-\tilde\beta_n^{\lambda_n}\right)\nabla h(x^*,\xi)\right]\\
&=&{\phi^*}''(0)E_P\left[\sqrt{\frac{\lambda_n}{{\phi^*}''(0)Var_P(h(x^*,\xi))}}(h(x^*,\xi)-E_P[h(x^*,\xi)])\nabla h(x^*,\xi)\right](1+o_p(1))\\
&=&\sqrt{\frac{\lambda_n{\phi^*}''(0)}{Var_P(h(x^*,\xi))}}Cov_P((h(x^*,\xi),\nabla h(x^*,\xi))(1+o_p(1))
\end{eqnarray*}
Putting this into \eqref{interim DRO proof4}, we get
\begin{eqnarray*}
&&-(E_P[\nabla^2h(x,\xi)])^{-1}\sqrt{\frac{\lambda_n{\phi^*}''(0)}{Var_P(h(x^*,\xi))}}Cov_P((h(x^*,\xi),\nabla h(x^*,\xi))-(E_P[\nabla^2h(x,\xi)])^{-1}E_{\hat P_n}[\nabla h(x^*,\xi)]{}\\
&&{}+o_p\left(\frac{1}{\sqrt n}+\sqrt\lambda_n\right)\\
&=&-(E_P[\nabla^2h(x^*,\xi)])^{-1}\langle\nabla h(x^*,\xi),\hat P_n-P\rangle-\sqrt{\lambda_n{\phi^*}''(0)}(E_P[\nabla^2h(x^*,\xi)])^{-1}\frac{Cov_P(h(x^*,\xi),\nabla h(x^*,\xi))}{\sqrt{Var_P(h(x^*,\xi))}}{}\\
&&{}+o_p\left(\frac{1}{\sqrt n}+\sqrt\lambda_n\right)
\end{eqnarray*}\hfill\Halmos

\endproof

\proof{Proof of Corollary \ref{Wasserstein thm}.}
The result is immediate by using, e.g., \cite{kuhn2019wasserstein} Theorem 10 to obtain \eqref{regularization Wasserstein}, on which we apply Theorem \ref{regularization thm}.
\endproof





\section{Conditions and Proof of Theorem \ref{thm:parametric}}\label{sec:parametric proof}
We make the following assumptions:
\begin{assumption}[Regularity conditions of parametric objective function]
We have:
\begin{enumerate}
\item The solution $x^*$ satisfies the first-order optimality condition, i.e., $\nabla_x\psi(x^*,\theta^*)=0$.
\item $\nabla_x\psi(\cdot,\cdot):\mathbb R^{d+m}\to\mathbb R^d$ is continuously differentiable at $(x^*,\theta^*)$ and the Hessian matrix of $\psi(\cdot,\cdot)$ with respect to $x$ at $(x^*,\theta^*)$,
$\nabla_x^2\psi(x^*,\theta^*)$, is invertible.
\end{enumerate}\label{assumption:parametric objective}
\end{assumption}

\begin{assumption}[Regularity conditions on penalty]
$R$ is twice continuously differentiable.\label{assumption:penalty}
\end{assumption}

We are now ready to prove Theorem \ref{thm:parametric}.

\proof{Proof of Theorem \ref{thm:parametric}.}
The result follows from a standard application of the delta method. By \eqref{parametric CLT}, we have $\hat\theta_n-\theta^*=O_p(1/\sqrt n)$ by the CLT and Slutsky's theorem. Now, by Assumption \ref{assumption:parametric objective}.1, $x^*$ is the solution to $\nabla_x\psi(x,\theta^*)+0\cdot\nabla_xR(x)=0$. By Assumptions \ref{assumption:parametric objective}.2 and \ref{assumption:penalty}, the function $\nabla_x\psi(x,\theta)+\lambda\nabla_xR(x)$ is continuously differentiable at $(x,\theta,\lambda)=(x^*,\theta^*,0)$, and $\nabla_x^2\psi(x^*,\theta^*)$ is invertible. Thus, by the implicit function theorem, we have a unique solution $x(\theta,\lambda)$ to the equation
$\nabla_x\psi(x,\theta)+\lambda\nabla_xg(x)=0$ for $(\theta,\lambda)$ in a small enough neighborhood of $(\theta,0)$, in which $x(\theta,\lambda)$ is continuously differentiable. Moreover, the gradient of $x(\theta,\lambda)$ at $(\theta^*,0)$ is given by $$(\nabla_x^2\psi(x^*,\theta^*))^{-1}[\nabla_{x\theta}\psi(x^*,\theta^*),\nabla_xR(x^*)]$$

Now let
$$r(h,\lambda)=x(\theta^*+h,\lambda)-x(\theta,0)-(\nabla_x^2\psi(x^*,\theta^*))^{-1}\nabla_{x\theta}\psi(x^*,\theta^*)h-(\nabla_x^2\psi(x^*,\theta^*))^{-1}\nabla_xR(x)\lambda$$
By the continuous differentiability of $x(\cdot,\cdot)$ argued above, we have
$r(h,\lambda)=o(\|(h,\lambda)\|)$ as $(h,\lambda)\to0$. Thus, by the definition of convergence in probability, we have
\begin{eqnarray*}
&&x(\hat\theta_n,\lambda)-x(\theta,0)-(\nabla_x^2\psi(x^*,\theta^*))^{-1}\nabla_{x\theta}\psi(x^*,\theta^*)(\hat\theta_n-\theta^*)-(\nabla_x^2\psi(x^*,\theta^*))^{-1}\nabla_xR(x^*)\lambda\\
&=&r(\hat\theta_n-\theta^*,\lambda)=o_p(\|(\hat\theta_n-\theta^*,\lambda)\|)=o_p(\|\hat\theta_n-\theta^*\|+|\lambda|)=o_p\left(\frac{1}{\sqrt n}+\lambda\right)
\end{eqnarray*}
Noting that $\hat x_n^{P-Reg,\lambda}=x(\hat\theta_n,\lambda)$, this concludes the theorem.\hfill\Halmos
\endproof

\section{Conditions and Proof of Theorem \ref{thm:Bayesian}}
\label{sec:Bayesian proof}
We make the following assumptions in addition to those in Section \ref{sec:parametric proof}:
\begin{assumption}[Additional regularity conditions for the objective function]
$\nabla_x\psi(x,\theta)$, $\nabla_{x\theta}\psi(x,\theta)$ and $\nabla_x^2\psi(x,\theta)$ are continuous and bounded uniformly over $x\in\mathcal X$ and $\theta\in\mathbb R^m$. Moreover, $\nabla_x^2\psi(x,\theta)$ is positive definite with eigenvalues bounded away from 0 and $\infty$ uniformly over $x\in\mathcal X$ and $\theta\in\mathbb R^m$.
\label{assumption:Bayesian}
\end{assumption}

\begin{assumption}[Additional regularity conditions for the penalty]
$\nabla_x^2R(x)$ is uniformly bounded over $x\in\mathcal X$.
\label{assumption:Bayesian penalty}
\end{assumption}

\begin{assumption}[Additional first-order optimality condition]
With probability converging to 1, there exists a solution 
$\hat x_n^{P-Bay,\lambda}$ as defined in \eqref{Bayesian solution} that satisfies the first-order optimality condition, i.e., it solves
$$\nabla_xE_{\Theta|D_n}[\psi(x,\Theta)]+\lambda\nabla_xR(x)$$for $n\to\infty$ and $\lambda\to0$.\label{assumption:Bayesian solution}
\end{assumption}

We are now ready to prove Theorem \ref{thm:Bayesian}.
\proof{Proof of Theorem \ref{thm:Bayesian}.}
By \eqref{posterior convergence}, we have
\begin{equation}
E_{\Theta|D_n}[\nabla_x\psi(x,\Theta)]=E_{\Theta|D_n}\left[\nabla_x\psi\left(x,\hat\theta_n+\frac{W_n}{\sqrt n}\right)\right]\label{Bayesian interim}
\end{equation}
where $\|W_n-N(0,J)\|_{TV}\stackrel{p}{\to}0$. Now consider \eqref{Bayesian interim} componentwise, i.e.,
\begin{equation}
E_{\Theta|D_n}\left[\nabla_x^j\psi\left(x,\hat\theta_n+\frac{W_n}{\sqrt n}\right)\right]\label{Bayesian interim1}
\end{equation}
where $\nabla_x^j$ for $j=1,\ldots,d$ denotes the $j$-th component of the gradient. By the mean value theorem and Assumption \ref{assumption:Bayesian}, we can write \eqref{Bayesian interim1} as
\begin{equation}
\nabla_x^j\psi(x,\hat\theta_n)+E_{\Theta|D_n}\left[\nabla_{x\theta}^j\psi(x,\hat\theta_n+\zeta_n^j)\frac{W_n}{\sqrt n}\right]\label{Bayesian interim2}
\end{equation}
where $\nabla_{x\theta}^j\psi(x,\theta)$ denotes the $j$-th row of $\nabla_{x\theta}\psi(x,\theta)$, and $\zeta_n^j$ is between 0 and $W_n/\sqrt n$. By the uniform boundedness of $\nabla_{x\theta}\psi(x,\theta)$ in Assumption \ref{assumption:Bayesian} and the almost sure uniform integrability of $W_n$, we have $\nabla_{x\theta}^j\psi(x,\hat\theta_n+\zeta_n^j)W_n$ also a.s. uniformly integrable. Thus, together with the continuity of $\nabla_{x\theta}\psi(x,\theta)$ in Assumption \ref{assumption:Bayesian}, $\hat\theta_n\stackrel{p}{\to}\theta^*$ as implied by \eqref{parametric CLT}, and $\|W_n-N(0,J)\|_{TV}\stackrel{p}{\to}0$ in \eqref{posterior convergence}, we have
$$E_{\Theta|D_n}[\nabla_{x\theta}^j\psi(x,\hat\theta_n+\zeta_n^j)W_n]\stackrel{p}{\to}0$$
Thus, 
\begin{equation}
E_{\Theta|D_n}[\nabla_x\psi(x,\Theta)]=\nabla_x\psi(x,\hat\theta_n)+o_p\left(\frac{1}{\sqrt n}\right)\label{Bayesian interim3}
\end{equation}
uniformly over $x\in\mathcal X$.

Now, define $\hat x_n^{P-Reg,\lambda}=x(\hat\theta,\lambda)$ as in the proof of Theorem \ref{thm:parametric}, and $\hat x^{P-Bay,\lambda}$ as the root of
$$E_{\Theta|D_n}[\nabla_x\psi(x,\Theta)]+\lambda\nabla_xR(x)$$
which exists with probability converging to 1 as $n\to\infty$ and coincides with the definition in Assumption \ref{assumption:Bayesian solution}, which is well-defined by using Assumption \ref{assumption:Bayesian}. 

Now, by the first-order optimality condition,
$$\nabla_x\psi(\hat x_n^{P-Reg,\lambda},\hat\theta_n)-\nabla_x\psi(\hat x_n^{P-EO},\hat\theta_n)=-\lambda\nabla_xR(\hat x_n^{P-Reg,\lambda})$$
By the delta method, and together with the uniform boundedness of $\nabla_x^2\psi(x,\theta)$ and the continuity of $\nabla_x^2\psi(x,\theta)$ in $x$ uniformly over $\theta$ in Assumption \ref{assumption:Bayesian}, we can rewrite the above as
$$\nabla_x^2\psi(\hat x_n^{P-EO},\hat\theta_n)(\hat x_n^{P-Reg,\lambda}-\hat x_n^{P-EO})+o_p(\|\hat x_n^{P-Reg,\lambda}-\hat x_n^{P-EO}\|)=-\lambda\nabla_xR(\hat x_n^{P-Reg,\lambda})$$
Note that then we have
\begin{align*}
\|\hat x_n^{P-Reg,\lambda}-\hat x_n^{P-EO}\|&\leq\|\nabla_x^2\psi(\hat x_n^{P-EO},\hat\theta_n)^{-1}\|\|\nabla_x^2\psi(\hat x_n^{P-EO},\hat\theta_n)(\hat x_n^{P-Reg,\lambda}-\hat x_n^{P-EO})\|\\
&=O(\lambda)+o_p(\|\hat x_n^{P-Reg,\lambda}-\hat x_n^{P-EO}\|)
\end{align*}
by using the uniform boundedness of the eigenvalues of $\nabla_x^2\psi(x,\theta)$ away from 0 and $\infty$ in Assumption \ref{assumption:Bayesian}, and the uniform boundedness of $\nabla_xR(x)$ in Assumption \ref{assumption:Bayesian penalty}. Thus $\hat x_n^{P-Reg,\lambda}-\hat x_n^{P-EO}\stackrel{p}{\to}0$ as $\lambda\to0$. We then have
\begin{equation}
\hat x_n^{P-Reg,\lambda}-\hat x_n^{P-EO}=-\lambda\nabla_x^2\psi(\hat x_n^{P-EO},\hat\theta_n)^{-1}\nabla_xR(\hat x_n^{P-Reg,\lambda})(1+o_p(1))\label{Bayesian interim4}
\end{equation}

Similarly, by the first-order optimality condition in Assumption \ref{assumption:Bayesian solution}, Assumption \ref{assumption:Bayesian} and \eqref{Bayesian interim3}, we have
$$\nabla_x\psi(\hat x_n^{P-Bay,\lambda},\hat\theta_n)-\nabla_x\psi(\hat x_n^{P-Bay,0},\hat\theta)=-\lambda\nabla_xR(\hat x_n^{P-Bay,\lambda})+o_p\left(\frac{1}{\sqrt n}\right)$$
So, by the delta method, and together with the uniform boundedness of $\nabla_x^2\psi(x,\theta)$ in Assumption \ref{assumption:Bayesian}, we have
$$\nabla_x^2\psi(\hat x_n^{P-Bay,0},\hat\theta_n)(\hat x_n^{P-Bay,\lambda}-\hat x_n^{P-Bay,0})+o_p(\|\hat x_n^{P-Bay,\lambda}-\hat x_n^{P-Bay,0}\|)=-\lambda\nabla_xR(\hat x_n^{P-Bay,\lambda})+o_p\left(\frac{1}{\sqrt n}\right)$$
Then we have
\begin{eqnarray*}
&&\|\hat x_n^{P-Bay,\lambda}-\hat x_n^{P-Bay,0}\|\leq\|\nabla_x^2\psi(\hat x_n^{P-Bay,0},\hat\theta_n)^{-1}\|\|\nabla_x^2\psi(\hat x_n^{P-Bay,0},\hat\theta_n)(\hat x_n^{P-Bay,\lambda}-\hat x_n^{P-Bay,0})\|\\
&=&O(\lambda)+o_p\left(\frac{1}{\sqrt n}\right)+o_p(\|\hat x_n^{P-Bay,\lambda}-\hat x_n^{P-Bay,0}\|)
\end{eqnarray*}
by using the uniform boundedness of the eigenvalues of $\nabla_x^2\psi(x,\theta)$ away from 0 and $\infty$ in Assumption \ref{assumption:Bayesian}, and the uniform boundedness of $\nabla_xR(x)$ in Assumption \ref{assumption:Bayesian penalty}. Thus $\hat x_n^{P-Bay,\lambda}-\hat x_n^{P-Bay,0}\stackrel{p}{\to}0$ as $\lambda\to0$ and $n\to\infty$. We then have
\begin{equation}
\hat x_n^{P-Bay,\lambda}-\hat x_n^{P-Bay,0}=-\lambda\nabla_x^2\psi(\hat x_n^{P-Bay,0},\hat\theta_n)^{-1}\nabla_xR(\hat x_n^{P-Bay,\lambda})+o_p\left(\frac{1}{\sqrt n}\right)\label{Bayesian interim5}
\end{equation}

Now, by the delta method like in the proof of Theorem \ref{thm:parametric}, and together with the continuity and uniform boundedness of $\nabla_x^2\psi(x,\theta)$ in Assumption \ref{assumption:Bayesian}, we have
$$\nabla_x\psi(\hat x_n^{P-Bay,0},\hat\theta_n)-\nabla_x\psi(\hat x_n^{P-EO},\hat\theta_n)=\nabla_x^2\psi(\hat x_n^{P-EO},\hat\theta_n)(\hat x_n^{P-Bay,0}-\hat x_n^{P-EO})+o_p(\|\hat x_n^{P-Bay,0}-\hat x_n^{P-EO}\|)$$
Note also that by \eqref{Bayesian interim3} and the first-order optimality conditions in Assumption \ref{assumption:parametric objective}.1 and Assumption \ref{assumption:Bayesian},
$$\nabla_x\psi(\hat x_n^{P-Bay,0},\hat\theta_n)+o_p\left(\frac{1}{\sqrt n}\right)=0=\nabla_x\psi(\hat x_n^{P-EO},\hat\theta_n)$$
so that
$$\nabla_x\psi(\hat x_n^{P-Bay,0},\hat\theta_n)-\nabla_x\psi(\hat x_n^{P-EO},\hat\theta_n)=o_p\left(\frac{1}{\sqrt n}\right)$$
Thus
\begin{align*}
\|\hat x_n^{P-Bay,0}-\hat x_n^{P-EO}\|&\leq\|\nabla_x^2\psi(\hat x_n^{P-EO},\hat\theta_n)^{-1}\|\|\nabla_x^2\psi(\hat x_n^{P-EO},\hat\theta_n)(\hat x_n^{P-Bay,0}-\hat x_n^{P-EO})\|\\
&=o_p\left(\frac{1}{\sqrt n}\right)+o_p(\|\hat x_n^{P-Bay,0}-\hat x_n^{P-EO}\|)
\end{align*}
by using the uniform boundedness of the eigenvalues of  $\nabla_x^2\psi(x,\theta)$ away from 0 and $\infty$ in Assumption \ref{assumption:Bayesian}. This gives
\begin{equation}
\hat x_n^{P-Bay,0}-\hat x_n^{P-EO}=o_p\left(\frac{1}{\sqrt n}\right)\label{Bayesian interim6}
\end{equation}

Putting together \eqref{Bayesian interim4}, \eqref{Bayesian interim5} and \eqref{Bayesian interim6}, we have
\begin{eqnarray*}
&&\hat x_n^{P-Bay,\lambda}-\hat x_n^{P-Reg,\lambda}\\
&=&(\hat x_n^{P-Bay,\lambda}-\hat x_n^{P-Bay,0})+(\hat x_n^{P-Bay,0}-\hat x_n^{P-EO})+(\hat x_n^{P-EO}-\hat x_n^{P-Reg,\lambda})\\
&=&-\lambda\nabla_x^2\psi(\hat x_n^{P-Bay,0},\hat\theta_n)^{-1}\nabla_xR(\hat x_n^{P-Bay,\lambda})+\lambda\nabla_x^2\psi(\hat x_n^{P-EO},\hat\theta_n)^{-1}\nabla_xR(\hat x_n^{P-Reg,\lambda})+o_p\left(\frac{1}{\sqrt n}\right)\\
&=&o_p\left(\frac{1}{\sqrt n}+\lambda\right)
\end{eqnarray*}
by the consistency of $\hat x_n^{P-Bay,\lambda},\hat x_n^{P-Reg,\lambda},\hat x_n^{P-Bay,0},\hat x_n^{P-EO}\stackrel{p}{\to} x^*$ implied by these equations and the continuity of $\nabla_x^2\psi(x,\theta)$ in Assumption \ref{assumption:Bayesian} and of $\nabla_xR(x)$ in Assumption \ref{assumption:Bayesian penalty}. Together with the conclusions from Theorem \ref{thm:parametric}, we prove the theorem.\hfill\Halmos
\endproof



\end{document}